\documentclass[final,3p]{elsarticle}
 \usepackage{graphics}
 \usepackage{graphicx}
 \usepackage{epsfig}
\usepackage{amssymb}
 \usepackage{amsthm}
 \usepackage{lineno}
 \usepackage{amsmath}
   \numberwithin{equation}{section}
\usepackage{mathrsfs}

\newtheorem{thm}{Theorem}[section]

\newtheorem{lem}[thm]{Lemma}

\newtheorem{defn}[thm]{Definition}

\biboptions{sort&compress,square}
\begin{document}
\begin{frontmatter}
\author{Tong Wu$^1$}
\ead{wut977@nenu.edu.cn}
\author{Jian Wang$^2$}
\ead{ wangj484@nenu.edu.cn}
\author{Yong Wang$^3$\corref{cor3}}
\ead{wangy581@nenu.edu.cn}
\cortext[cor3]{Corresponding author.}

\address{1.School of Mathematics and Statistics, Northeast Normal University,
Changchun, 130024, China}
\address{2.School of Science, Tianjin University of Technology and Education,
Tianjin, 300222, China}
\address{3.School of Mathematics and Statistics, Northeast Normal University,
Changchun, 130024, China}

\title{ Dirac-Witten Operators and the Kastler-Kalau-Walze type theorem for manifolds with boundary }
\begin{abstract}
In this paper, we obtain two Lichnerowicz type formulas for the Dirac-Witten operators. And we give the proof of Kastler-Kalau-Walze type theorems for the Dirac-Witten operators on 4-dimensional and 6-dimensional compact manifolds with (resp.without) boundary.
\end{abstract}
\end{frontmatter}
\section{Introduction}
 Until now, many geometers have studied noncommutative residues. In \cite{Gu,Wo}, authors found noncommutative residues are of great importance to the study of noncommutative geometry. In \cite{Co1}, Connes used the noncommutative residue to derive a conformal 4-dimensional Polyakov action analogy. Connes showed us that the noncommutative residue on a compact manifold $M$ coincided with the Dixmier's trace on pseudodifferential operators of order $-{\rm {dim}}M$ in \cite{Co2}.
And Connes claimed the noncommutative residue of the square of the inverse of the Dirac operator was proportioned to the Einstein-Hilbert action.  Kastler \cite{Ka} gave a
brute-force proof of this theorem. Kalau and Walze proved this theorem in the normal coordinates system simultaneously in \cite{KW} .
Ackermann proved that
the Wodzicki residue  of the square of the inverse of the Dirac operator ${\rm  Wres}(D^{-2})$ in turn is essentially the second coefficient
of the heat kernel expansion of $D^{2}$ in \cite{Ac}.

On the other hand, Wang generalized the Connes' results to the case of manifolds with boundary in \cite{Wa1,Wa2},
and proved the Kastler-Kalau-Walze type theorem for the Dirac operator and the signature operator on lower-dimensional manifolds
with boundary \cite{Wa3}. In \cite{Wa3,Wa4}, Wang computed $\widetilde{{\rm Wres}}[\pi^+D^{-1}\circ\pi^+D^{-1}]$ and $\widetilde{{\rm Wres}}[\pi^+D^{-2}\circ\pi^+D^{-2}]$, where the two operators are symmetric, in these cases the boundary term vanished. But for $\widetilde{{\rm Wres}}[\pi^+D^{-1}\circ\pi^+D^{-3}]$, Wang got a nonvanishing boundary term \cite{Wa5}, and give a theoretical explanation for gravitational action on boundary. In others words, Wang provides a kind of method to study the Kastler-Kalau-Walze type theorem for manifolds with boundary.
In \cite{lkl}, L\'{o}pez and his collaborators introduced an elliptic differential operator which is called the Novikov operator. In \cite{WW}, Wei and Wang proved Kastler-Kalau-Walze type theorem for modified Novikov operators on compact manifolds. In \cite{zx}, in order to prove the nonsymmetric positive mass theorem, Zhang introduced the Dirac-Witten operator. The motivation of this paper is
to prove the Kastler-Kalau-Walze type theorem for the Dirac-Witten operators. \\
\indent The paper is organized in the following way. In Section 2, by using the definition of Dirac-Witten operators, we compute the Lichnerowicz formulas for the Dirac-Witten operators. In Section 3 and in Section 4,
 we prove the Kastler-Kalau-Walze type theorem for 4-dimensional and 6-dimensional manifolds with boundary for the Dirac-Witten operators respectively.
\section{The Dirac-Witten Operators and their Lichnerowicz formulas}

Firstly we introduce some notations about the Dirac-Witten Operators. Let $M$ be a $n$-dimensional ($n\geq 3$) oriented compact spin Riemannian manifold with a Riemannian metric $g^{M}$. And let $\nabla^L$ be the Levi-Civita connection about $g^M$. In the local coordinates $\{x_i; 1\leq i\leq n\}$ and the
fixed orthonormal frame $\{e_1,\cdots,e_n\}$, the connection matrix $(\omega_{s,t})$ is defined by
\begin{equation}
\nabla^L(e_1,\cdots,e_n)= (e_1,\cdots,e_n)(\omega_{s,t}).
\end{equation}
\indent Let $c(e_j)$ be the Clifford action.
Suppose that $\partial_{i}$ is a natural local frame on $TM$
and $(g^{ij})_{1\leq i,j\leq n}$ is the inverse matrix associated to the metric
matrix  $(g_{ij})_{1\leq i,j\leq n}$ on $M$.
By \cite{Wa3}, we have the Dirac operator
\begin{align}
D&=\sum^n_{i=1}c(e_i)\bigg[e_i-\frac{1}{4}\sum_{s,t}\omega_{s,t}
(e_i)c(e_s)c(e_t)\bigg]
\end{align}

\indent Then the Dirac-Witten operators $\widetilde{D}$ and ${\widetilde{D}}^*$ are defined by
\begin{align}
\widetilde{D}&=D+f_1\sum_{u<v}(p_{uv}-p_{vu})c(e_u)c(e_v)+f_2\\
&=\sum^n_{i=1}c(e_i)\bigg[e_i-\frac{1}{4}\sum_{s,t}\omega_{s,t}
(e_i)c(e_s)c(e_t)\bigg]+f_1\sum_{u<v}(p_{uv}-p_{vu})c(e_u)c(e_v)+f_2,\notag\\
{\widetilde{D}}^*&=D-\overline{f_1}\sum_{u<v}(p_{uv}-p_{vu})c(e_u)c(e_v)+\overline{f_2}\nonumber\\
&=\sum^n_{i=1}c(e_i)\bigg[e_i-\frac{1}{4}\sum_{s,t}\omega_{s,t}
(e_i)c(e_s)c(e_t)\bigg]-\overline{f_1}\sum_{u<v}(p_{uv}-p_{vu})c(e_u)c(e_v)+\overline{f_2}.\notag
\end{align}
where $f_1,f_2$ is a complex number and $p_{uv}=p(e_u,e_v)$, $p$ is a $(0,2)$-tensor. Then when $f_1=\frac{\sqrt{-1}}{2}$, $f_2=-\frac{\sqrt{-1}}{2}\sum_{i}p_{ii}$, $\widetilde{D}$ is the Dirac-Witten Operator defined by \cite{zx}.\\
\indent Next by computations, we get Lichnerowicz formulas,
\begin{thm} The following equalities hold:
\begin{align}
{\widetilde{D}}^*\widetilde{D}&=-\Big[g^{ij}(\nabla_{\partial_{i}}\nabla_{\partial_{j}}-
\nabla_{\nabla^{L}_{\partial_{i}}\partial_{j}})\Big]+\frac{1}{4}s+(f_1\overline{f_2}-\overline{f_1}f_2)\sum_{u<v}(p_{uv}-p_{vu})c(e_u)c(e_v)\nonumber\\
&+\frac{1}{4}\sum_{i}[c(e_{i})(f_1\sum_{u<v}(p_{uv}-p_{vu})c(e_u)c(e_v)+f_2)+(-\overline{f_1}\sum_{u<v}(p_{uv}-p_{vu})c(e_u)c(e_v)+\overline{f_2})c(e_i)]^2\nonumber\\
&+\frac{1}{2}[c(e_{j})e_j(f_1\sum_{u<v}(p_{uv}-p_{vu})c(e_u)c(e_v)+f_2)-e_j(-\overline{f_1}\sum_{u<v}(p_{uv}-p_{vu})c(e_u)c(e_v)+\overline{f_2})c(e_j)]\nonumber\\
&-f_1\overline{f_1}[\sum_{u<v}(p_{uv}-p_{vu})c(e_u)c(e_v)]^2+f_2\overline{f_2},\nonumber\\
{\widetilde{D}}^2&=-\Big[g^{ij}(\nabla_{\partial_{i}}\nabla_{\partial_{j}}-
\nabla_{\nabla^{L}_{\partial_{i}}\partial_{j}})\Big]+\frac{1}{4}s+{f_1}^2[\sum_{u<v}(p_{uv}-p_{vu})c(e_u)c(e_v)]^2\nonumber\\
&+\frac{1}{4}\sum_{i}[c(e_{i})(f_1\sum_{u<v}(p_{uv}-p_{vu})c(e_u)c(e_v)+f_2)+(f_1\sum_{u<v}(p_{uv}-p_{vu})c(e_u)c(e_v)+f_2)c(e_i)]^2\nonumber\\
&-\frac{1}{2}[e_j(f_1\sum_{u<v}(p_{uv}-p_{vu})c(e_u)c(e_v)+f_2)c(e_{j})-c(e_{j})e_j(f_1\sum_{u<v}(p_{uv}-p_{vu})c(e_u)c(e_v)+f_2)]\nonumber\\
&+2f_1f_2\sum_{u<v}(p_{uv}-p_{vu})c(e_u)c(e_v)+{f_2}^2.\nonumber\\
\end{align}
where $s$ is the scalar curvature.
\end{thm}
Let $M$ be a smooth compact oriented spin Riemannian $n$-dimensional manifolds without boundary and $N$ be a vector bundle on $M$. If $P$ is a differential operator of Laplace type, then it has locally the form
\begin{equation}
P=-(g^{ij}\partial_i\partial_j+A^i\partial_i+B),
\end{equation}
where $\partial_{i}$  is a natural local frame on $TM$
and $(g^{ij})_{1\leq i,j\leq n}$ is the inverse matrix associated to the metric
matrix  $(g_{ij})_{1\leq i,j\leq n}$ on $M$,
 and $A^{i}$ and $B$ are smooth sections
of $\textrm{End}(N)$ on $M$ (endomorphism). If a Laplace type
operator $P$  satisfies (2.7), then there is a unique
connection $\nabla$ on $N$ and a unique endomorphism $E$ such that
 \begin{equation}
P=-\Big[g^{ij}(\nabla_{\partial_{i}}\nabla_{\partial_{j}}-
 \nabla_{\nabla^{L}_{\partial_{i}}\partial_{j}})+E\Big],
\end{equation}
where $\nabla^{L}$ is the Levi-Civita connection on $M$. Moreover
(with local frames of $T^{*}M$ and $N$), $\nabla_{\partial_{i}}=\partial_{i}+\omega_{i} $
and $E$ are related to $g^{ij}$, $A^{i}$ and $B$ through
 \begin{eqnarray}
&&\omega_{i}=\frac{1}{2}g_{ij}\big(A^{i}+g^{kl}\Gamma_{ kl}^{j} \texttt{id}\big),\\
&&E=B-g^{ij}\big(\partial_{i}(\omega_{j})+\omega_{i}\omega_{j}-\omega_{k}\Gamma_{ ij}^{k} \big),
\end{eqnarray}
where $\Gamma_{ kl}^{j}$ is the  Christoffel coefficient of $\nabla^{L}$.\\
\indent Let $g^{ij}=g(dx_{i},dx_{j})$, $\xi=\sum_{k}\xi_{j}dx_{j}$ and $\nabla^L_{\partial_{i}}\partial_{j}=\sum_{k}\Gamma_{ij}^{k}\partial_{k}$,  we denote that
\begin{align}
&\sigma_{i}=-\frac{1}{4}\sum_{s,t}\omega_{s,t}
(e_i)c(e_s)c(e_t)
;\nonumber\\
&\xi^{j}=g^{ij}\xi_{i};~~~~\Gamma^{k}=g^{ij}\Gamma_{ij}^{k};~~~~\sigma^{j}=g^{ij}\sigma_{i}.
\end{align}
\indent Then the Dirac-Witten operators $\widetilde{D}$ and ${\widetilde{D}}^*$ can be written as
\begin{align}
\widetilde{D}=\sum^n_{i=1}c(e_i)[e_i+\sigma_{i}]
+f_1\sum_{u<v}(p_{uv}-p_{vu})c(e_u)c(e_v)+f_2;\nonumber\\
{\widetilde{D}}^*=\sum^n_{i=1}c(e_i)[e_i+\sigma_{i}]
-\overline{f_1}\sum_{u<v}(p_{uv}-p_{vu})c(e_u)c(e_v)+\overline{f_2}.\nonumber\\
\end{align}
\indent By \cite{Ka}, we have
\begin{align}
D^2&=-\Delta_0+\frac{1}{4}s\\\notag
&=-g^{ij}(\nabla_i^L\nabla_j^L-\Gamma_{ij}^k\nabla_k^L)+\frac{1}{4}s\\\notag
&=-\sum_{ij}g^{ij}[\partial_i\partial_j+2\sigma_i\partial_j-\Gamma_{ij}^k\partial_k+\partial_i\sigma_j+\sigma_i\sigma_j-\Gamma_{ij}^k\sigma_k]+\frac{1}{4}s.\\\notag
\end{align}
\indent By (2.10), we have
\begin{align}
{\widetilde{D}}^*D&=D^2+D[f_1\sum_{u<v}(p_{uv}-p_{vu})c(e_u)c(e_v)+f_2]+[-\overline{f_1}\sum_{u<v}(p_{uv}-p_{vu})\nonumber\\
&c(e_u)c(e_v)+\overline{f_2}]D+[f_1\sum_{u<v}(p_{uv}-p_{vu})c(e_u)c(e_v)+f_2]\nonumber\\
&[-\overline{f_1}\sum_{u<v}(p_{uv}-p_{vu})c(e_u)c(e_v)+\overline{f_2}],\nonumber\\
\end{align}
\begin{align}
&D[f_1\sum_{u<v}(p_{uv}-p_{vu})c(e_u)c(e_v)+f_2]+[-\overline{f_1}\sum_{u<v}(p_{uv}-p_{vu})c(e_u)c(e_v)+\overline{f_2}]D\nonumber\\
&=\sum_{i,j}g^{i,j}\Big[c(\partial_{i})(f_1\sum_{u<v}(p_{uv}-p_{vu})c(e_u)c(e_v)+f_2)
+(-\overline{f_1}\sum_{u<v}(p_{uv}-p_{vu})c(e_u)c(e_v)+\overline{f_2})c(\partial_{i})\Big]\partial_{j}\nonumber\\
&-\sum_{i,j}g^{i,j}\Big[(\overline{f_1}\sum_{u<v}(p_{uv}-p_{vu})c(e_u)c(e_v)-\overline{f_2})c(\partial_{i})\sigma_{j}-c(\partial_{i})\partial_{j}(f_1\sum_{u<v}(p_{uv}-p_{vu})c(e_u)c(e_v)+f_2)\nonumber\\
&-c(\partial_{i})\sigma_j(f_1\sum_{u<v}(p_{uv}-p_{vu})c(e_u)c(e_v)+f_2)\Big],\nonumber\\
\end{align}
then we obtain
\begin{align}
{\widetilde{D}}^*\widetilde{D}&=-\sum_{i,j}g^{i,j}\Big[\partial_{i}\partial_{j}+2\sigma_{i}\partial_{j}-\Gamma_{i,j}^{k}\partial_{k}+\partial_{i}\sigma_{j}
+\sigma_{i}\sigma_{j} -\Gamma_{i,j}^{k}\sigma_{k}\Big]+\frac{1}{4}s+\sum_{i,j}g^{i,j}\Big[c(\partial_{i})\nonumber\\
&(f_1\sum_{u<v}(p_{uv}-p_{vu})c(e_u)c(e_v)+f_2)
+(-\overline{f_1}\sum_{u<v}(p_{uv}-p_{vu})c(e_u)c(e_v)+\overline{f_2})c(\partial_{i})\Big]\partial_{j}\nonumber\\
&-\sum_{i,j}g^{i,j}\Big[(\overline{f_1}\sum_{u<v}(p_{uv}-p_{vu})c(e_u)c(e_v)-\overline{f_2})c(\partial_i)\sigma_j-c(\partial_i)\partial_j(f_1\sum_{u<v}(p_{uv}-p_{vu})\nonumber\\
&c(e_u)c(e_v)+f_2)-c(\partial_{i})\sigma_{j}(f_1\sum_{u<v}(p_{uv}-p_{vu})c(e_u)c(e_v)+f_2)\Big]-f_1\overline{f_1}[\sum_{u<v}(p_{uv}-p_{vu})\nonumber\\
&c(e_u)c(e_v)]^2+(f_1\overline{f_2}-\overline{f_1}f_2)\sum_{u<v}(p_{uv}-p_{vu})c(e_u)c(e_v)+f_2\overline{f_2}.\nonumber\\
\end{align}
\indent Similarly, we have
\begin{align}
{\widetilde{D}}^2&=-\sum_{i,j}g^{i,j}\Big[\partial_{i}\partial_{j}+2\sigma_{i}\partial_{j}-\Gamma_{i,j}^{k}\partial_{k}+\partial_{i}\sigma_{j}
+\sigma_{i}\sigma_{j}-\Gamma_{i,j}^{k}\sigma_{k}\Big]+\frac{1}{4}s+\sum_{i,j}g^{i,j}\Big[c(\partial_{i})\nonumber\\
&(f_1\sum_{u<v}(p_{uv}-p_{vu})c(e_u)c(e_v)+f_2)+(f_1\sum_{u<v}(p_{uv}-p_{vu})c(e_u)c(e_v)+f_2)c(\partial_{i})\Big]\partial_{j}\nonumber\\
&+\sum_{i,j}g^{i,j}\Big[(f_1\sum_{u<v}(p_{uv}-p_{vu})c(e_u)c(e_v)+f_2)c(\partial_{i})\sigma_{j}+c(\partial_{i})\partial_{j}(f_1\sum_{u<v}(p_{uv}-p_{vu})\nonumber\\
&c(e_u)c(e_v)+f_2)+c(\partial_{i})\sigma_{j}(f_1\sum_{u<v}(p_{uv}-p_{vu})c(e_u)c(e_v)+f_2)\Big]+f_1^2[\sum_{u<v}(p_{uv}-p_{vu})\nonumber\\
&c(e_u)c(e_v)]^2+2f_1f_2\sum_{u<v}(p_{uv}-p_{vu})c(e_u)c(e_v)+f_2^2.\nonumber\\
\end{align}
\indent By (2.6), (2.7), (2.8) and (2.14), we have
\begin{align}
(\omega_{i})_{{\widetilde{D}}^*{\widetilde{D}}}&=\sigma_{i}-\frac{1}{2}\Big[c(\partial_{i})(f_1\sum_{u<v}(p_{uv}-p_{vu})c(e_u)c(e_v)+f_2)\nonumber\\
&+(-\overline{f_1}\sum_{u<v}(p_{uv}-p_{vu})c(e_u)c(e_v)+\overline{f_2})c(\partial_{i})\Big].\nonumber\\
\end{align}
\begin{align}
E_{{\widetilde{D}}^*\widetilde{D}}&=-c(\partial_{i})\sigma^{i}(f_1\sum_{u<v}(p_{uv}-p_{vu})c(e_u)c(e_v)+f_2)-\frac{1}{4}s+(\overline{f_1}\sum_{u<v}(p_{uv}-p_{vu})c(e_u)c(e_v)-\overline{f_2})\nonumber\\
&c(\partial_i)\sigma^i+c(\partial_i)\partial^i(f_1\sum_{u<v}(p_{uv}-p_{vu})c(e_u)c(e_v)+f_2)+\frac{1}{2}\partial^j[c(\partial_j)(f_1\sum_{u<v}(p_{uv}-p_{vu})c(e_u)c(e_v)\nonumber\\
&+f_2)+(-\overline{f_1}\sum_{u<v}(p_{uv}-p_{vu})c(e_u)c(e_v)+\overline{f_2})c(\partial_j)]-\frac{1}{2}[c(\partial_j)(f_1\sum_{u<v}(p_{uv}-p_{vu})c(e_u)c(e_v)+f_2)\nonumber\\
&+(-\overline{f_1}\sum_{u<v}(p_{uv}-p_{vu})c(e_u)c(e_v)+\overline{f_2})c(\partial_j)]\sigma^j-\frac{g^{ij}}{4}[c(\partial_{i})(f_1\sum_{u<v}(p_{uv}-p_{vu})c(e_u)c(e_v)+f_2)\nonumber\\
&+
(-\overline{f_1}\sum_{u<v}(p_{uv}-p_{vu})c(e_u)c(e_v)+\overline{f_2})c(\partial^i)]\cdot[c(\partial_{j})(f_1\sum_{u<v}(p_{uv}-p_{vu})c(e_u)c(e_v)+f_2)\nonumber\\
&-(\overline{f_1}\sum_{u<v}(p_{uv}-p_{vu})c(e_u)c(e_v)-\overline{f_2})c(\partial_{j})]-\frac{1}{2}\Gamma^{k}[c(\partial_{k})(f_1\sum_{u<v}(p_{uv}-p_{vu})c(e_u)c(e_v)+f_2)\nonumber\\
&+(-\overline{f_1}\sum_{u<v}(p_{uv}-p_{vu})c(e_u)c(e_v)+\overline{f_2})c(\partial_{k})]-\frac{1}{2}\sigma^{j}[c(\partial_{j})(f_1\sum_{u<v}(p_{uv}-p_{vu})c(e_u)c(e_v)+f_2)\nonumber\\
&+(-\overline{f_1}\sum_{u<v}(p_{uv}-p_{vu})c(e_u)c(e_v)+\overline{f_2})c(\partial_{j})]+f_1\overline{f_1}[\sum_{u<v}(p_{uv}-p_{vu})c(e_u)c(e_v)]^2-(f_1\overline{f_2}-\overline{f_1}f_2)\nonumber\\
&\sum_{u<v}(p_{uv}-p_{vu})c(e_u)c(e_v)-f_2\overline{f_2}.\nonumber\\
\end{align}
\indent Since $E$ is globally
defined on $M$, taking normal coordinates at $x_0$, we have
$\sigma^{i}(x_0)=0$, $\partial^{j}[c(\partial_{j})](x_0)=0$,
$\Gamma^k(x_0)=0$, $g^{ij}(x_0)=\delta^j_i$, then
\begin{align}
E_{{\widetilde{D}}^*\widetilde{D}}(x_0)&=-\frac{1}{4}s-(f_1\overline{f_2}-\overline{f_1}f_2)\sum_{u<v}(p_{uv}-p_{vu})c(e_u)c(e_v)+f_1\overline{f_1}[\sum_{u<v}(p_{uv}-p_{vu})c(e_u)c(e_v)]^2\nonumber\\
&-f_2\overline{f_2}-\frac{1}{4}\sum_{i}[c(e_{i})(f_1\sum_{u<v}(p_{uv}-p_{vu})c(e_u)c(e_v)+f_2)+(-\overline{f_1}\sum_{u<v}(p_{uv}-p_{vu})\nonumber\\
&c(e_u)c(e_v)+\overline{f_2})c(e_i)]^2-\frac{1}{2}[c(e_{j})\nabla^{\bigwedge^*T^*M}_{e_j}(f_1\sum_{u<v}(p_{uv}-p_{vu})c(e_u)c(e_v)+f_2)\nonumber\\
&-\nabla^{\bigwedge^*T^*M}_{e_j}(-\overline{f_1}\sum_{u<v}(p_{uv}-p_{vu})c(e_u)c(e_v)+\overline{f_2})c(e_j)].\nonumber\\
\end{align}
\indent Similarly, we have
\begin{align}
E_{{\widetilde{D}}^2}(x_0)&=-\frac{1}{4}s-f_1^2[\sum_{u<v}(p_{uv}-p_{vu})c(e_u)c(e_v)]^2-2f_1f_2\sum_{u<v}(p_{uv}-p_{vu})c(e_u)c(e_v)\nonumber\\
&-f_2^2-\frac{1}{4}\sum_{i}[c(e_{i})(f_1\sum_{u<v}(p_{uv}-p_{vu})c(e_u)c(e_v)+f_2)+(f_1\sum_{u<v}(p_{uv}-p_{vu})\nonumber\\
&c(e_u)c(e_v)+f_2)c(e_i)]^2+\frac{1}{2}[\nabla^{\bigwedge^*T^*M}_{e_j}(f_1\sum_{u<v}(p_{uv}-p_{vu})c(e_u)c(e_v)+f_2)c(e_j)\nonumber\\
&-c(e_{j})\nabla^{\bigwedge^*T^*M}_{e_j}(f_1\sum_{u<v}(p_{uv}-p_{vu})c(e_u)c(e_v)+f_2)],\nonumber\\
\end{align}
by (2.5),we get Theorem 2.1.\\
\indent From \cite{Ac}, we konw that the noncommutative residue of a generalized laplacian $\overline{\Delta}$ is expressed as
\begin{equation}
(n-2)\Phi_{2}(\overline{\Delta})=(4\pi)^{-\frac{n}{2}}\Gamma(\frac{n}{2})\widetilde{res}(\overline{\Delta}^{-\frac{n}{2}+1}),
\end{equation}
where $\Phi_{2}(\overline{\Delta})$ denotes the integral over the diagonal part of the second
coefficient of the heat kernel expansion of $\overline{\Delta}$.
Now let $\overline{\Delta}={\widetilde{D}}^*{\widetilde{D}}$ and ${\widetilde{D}}^*{\widetilde{D}}=\Delta-E$, then we have
\begin{align}
{\rm Wres}({\widetilde{D}}^*{\widetilde{D}})^{-\frac{n-2}{2}}
=\frac{(n-2)(4\pi)^{\frac{n}{2}}}{(\frac{n}{2}-1)!}\int_{M}{\rm tr}(\frac{1}{6}s+E_{\widetilde{D}^*\widetilde{D}})d{\rm Vol_{M} },
\end{align}
\begin{align}
{\rm Wres}({\widetilde{D}}^2)^{-\frac{n-2}{2}}
=\frac{(n-2)(4\pi)^{\frac{n}{2}}}{(\frac{n}{2}-1)!}\int_{M}{\rm tr}(\frac{1}{6}s+E_{{\widetilde{D}}^2})d{\rm Vol_{M} },
\end{align}
where ${\rm Wres}$ denote the noncommutative residue.
By computations, we have
\begin{align}
&tr[\sum_{u<v}(p_{uv}-p_{vu})c(e_u)c(e_v)]=0;\nonumber\\
&tr[\sum_{u<v}(p_{uv}-p_{vu})c(e_u)c(e_v)]^2=-\sum_{u<v}(p_{uv}-p_{vu})^2{\rm tr}[{\rm \texttt{id}}];\nonumber\\
&tr\bigg(\sum_{i}[c(e_{i})(f_1\sum_{u<v}(p_{uv}-p_{vu})c(e_u)c(e_v)+f_2)]^2\bigg)\nonumber\\
&=[f_1^2(n-4)\sum_{u<v}(p_{uv}-p_{vu})^2-nf_2^2]{\rm tr}[{\rm \texttt{id}}];\nonumber\\
&tr\bigg(c(e_{j})\nabla^{\bigwedge^*T^*M}_{e_j}(f_1\sum_{u<v}(p_{uv}-p_{vu})c(e_u)c(e_v)+f_2)\bigg)=0.\nonumber\\
\end{align}
Then by (2.23), we get
\begin{align}
tr(E_{{\widetilde{D}}^*\widetilde{D}})&=\bigg[-\frac{s}{4}-\frac{1}{4}\{[(f_1^2+\overline{f_1}^2)(n-4)-2nf_1\overline{f_1}]\sum_{u<v}(p_{uv}-p_{vu})^2\nonumber\\
&+2nf_2\overline{f_2}-nf_2^2-n\overline{f_2}\}-f_2\overline{f_2}-f_1\overline{f_1}\sum_{u<v}(p_{uv}-p_{vu})^2\bigg]{\rm tr}[{\rm \texttt{id}}],\nonumber\\
tr(E_{{\widetilde{D}}^2})&=[-\frac{s}{4}+(3-n)f_1^2\sum_{u<v}(p_{uv}-p_{vu})^2+(n-1)f_2^2]{\rm tr}[{\rm \texttt{id}}].\nonumber\\
\end{align}
\begin{thm} If $M$ is a $n$-dimensional compact oriented spin manifolds without boundary, and $n$ is even, then we get the following equalities :
\begin{align}
{\rm Wres}({\widetilde{D}}^*\widetilde{D})^{-\frac{n-2}{2}}
&=\frac{(n-2)(4\pi)^{\frac{n}{2}}}{(\frac{n}{2}-1)!}\int_{M}2^\frac{n}{2}\bigg(
-\frac{1}{12}s-\frac{1}{4}\{[(f_1^2+\overline{f_1}^2)(n-4)-2nf_1\overline{f_1}]\nonumber\\
&\sum_{u<v}(p_{uv}-p_{vu})^2+2nf_2\overline{f_2}-nf_2^2-n\overline{f_2}^2\}-f_2\overline{f_2}-f_1\overline{f_1}\sum_{u<v}(p_{uv}-p_{vu})^2\bigg)d{\rm Vol_{M}}.\nonumber\\
\end{align}
\begin{align}
{\rm Wres}({\widetilde{D}}^2)^{-\frac{n-2}{2}}
&=\frac{(n-2)(4\pi)^{\frac{n}{2}}}{(\frac{n}{2}-1)!}\int_{M}2^\frac{n}{2}\bigg(
-\frac{1}{12}s+(3-n)f_1^2\sum_{u<v}(p_{uv}-p_{vu})^2+(n-1)f_2^2\bigg)d{\rm Vol_{M}}.\nonumber\\
\end{align}
where $s$ is the scalar curvature.
\end{thm}

\section{A Kastler-Kalau-Walze type theorem for $4$-dimensional manifolds with boundary}
 We firstly recall that some basic facts and formulas about Boutet de
Monvel's calculus and the definition of the noncommutative residue for manifolds with boundary which will be used in the following. For more details, (see in Section 2 in \cite{Wa3}).\\
 \indent Let $U\subset M$ be a collar neighborhood of $\partial M$ which is diffeomorphic with $\partial M\times [0,1)$. By the definition of $h(x_n)\in C^{\infty}([0,1))$
and $h(x_n)>0$, there exists $\widehat{h}\in C^{\infty}((-\varepsilon,1))$ such that $\widehat{h}|_{[0,1)}=h$ and $\widehat{h}>0$ for some
sufficiently small $\varepsilon>0$. Then there exists a metric $g'$ on $\widetilde{M}=M\bigcup_{\partial M}\partial M\times
(-\varepsilon,0]$ which has the form on $U\bigcup_{\partial M}\partial M\times (-\varepsilon,0 ]$
\begin{equation}
g'=\frac{1}{\widehat{h}(x_{n})}g^{\partial M}+dx _{n}^{2} ,
\end{equation}
such that $g'|_{M}=g$. We fix a metric $g'$ on the $\widetilde{M}$ such that $g'|_{M}=g$.

Let Fourier transformation $F'$  be
\begin{equation}
F':L^2({\bf R}_t)\rightarrow L^2({\bf R}_v);~F'(u)(v)=\int e^{-ivt}u(t)dt\\
\end{equation}
and let
\begin{equation}
r^{+}:C^\infty ({\bf R})\rightarrow C^\infty (\widetilde{{\bf R}^+});~ f\rightarrow f|\widetilde{{\bf R}^+};~
\widetilde{{\bf R}^+}=\{x\geq0;x\in {\bf R}\}.
\end{equation}
 where $\Phi({\bf R})$
denotes the Schwartz space and $\Phi(\widetilde{{\bf R}^+}) =r^+\Phi({\bf R})$, $\Phi(\widetilde{{\bf R}^-}) =r^-\Phi({\bf R})$.\\
\indent We define $H^+=F'(\Phi(\widetilde{{\bf R}^+}));~ H^-_0=F'(\Phi(\widetilde{{\bf R}^-}))$ which satisfies
$H^+\bot H^-_0$. We have the following
 property: $h\in H^+~(H^-_0)$ if and only if $h\in C^\infty({\bf R})$ which has an analytic extension to the lower (upper) complex
half-plane $\{{\rm Im}\xi<0\}~(\{{\rm Im}\xi>0\})$ such that for all nonnegative integer $l$,
 \begin{equation}
\frac{d^{l}h}{d\xi^l}(\xi)\sim\sum^{\infty}_{k=1}\frac{d^l}{d\xi^l}(\frac{c_k}{\xi^k}),
\end{equation}
as $|\xi|\rightarrow +\infty,{\rm Im}\xi\leq0~({\rm Im}\xi\geq0)$.
Let $H'$ be the space of all polynomials and $H^-=H^-_0\bigoplus H';~H=H^+\bigoplus H^-.$ Denote by $\pi^+~(\pi^-)$ respectively the
 projection on $H^+~(H^-)$. For calculations, we take $H=\widetilde H=\{$rational functions having no poles on the real axis$\}$ ($\tilde{H}$
 is a dense set in the topology of $H$). Then on $\tilde{H}$,
 \begin{equation}
\pi^+h(\xi_0)=\frac{1}{2\pi i}\lim_{u\rightarrow 0^{-}}\int_{\Gamma^+}\frac{h(\xi)}{\xi_0+iu-\xi}d\xi,
\end{equation}
where $\Gamma^+$ is a Jordan close curve
included ${\rm Im}(\xi)>0$ surrounding all the singularities of $h$ in the upper half-plane and
$\xi_0\in {\bf R}$. Similarly, define $\pi'$ on $\tilde{H}$,
\begin{equation}
\pi'h=\frac{1}{2\pi}\int_{\Gamma^+}h(\xi)d\xi.
\end{equation}
So, $\pi'(H^-)=0$. For $h\in H\bigcap L^1({\bf R})$, $\pi'h=\frac{1}{2\pi}\int_{{\bf R}}h(v)dv$ and for $h\in H^+\bigcap L^1({\bf R})$, $\pi'h=0$.

Let $M$ be a $n$-dimensional compact oriented spin manifold with boundary $\partial M$.
Denote by $\mathcal{B}$ Boutet de Monvel's algebra, we recall the main theorem in \cite{FGLS,Wa3}.
\begin{thm}\label{th:32}{\rm\cite{FGLS}}{\bf(Fedosov-Golse-Leichtnam-Schrohe)}
 Let $X$ and $\partial X$ be connected, ${\rm dim}X=n\geq3$,
 $A=\left(\begin{array}{lcr}\pi^+P+G &   K \\
T &  S    \end{array}\right)$ $\in \mathcal{B}$ , and denote by $p$, $b$ and $s$ the local symbols of $P,G$ and $S$ respectively.
 Define:
 \begin{align}
{\rm{\widetilde{Wres}}}(A)&=\int_X\int_{\bf S}{\rm{tr}}_E\left[p_{-n}(x,\xi)\right]\sigma(\xi)dx \nonumber\\
&+2\pi\int_ {\partial X}\int_{\bf S'}\left\{{\rm tr}_E\left[({\rm{tr}}b_{-n})(x',\xi')\right]+{\rm{tr}}
_F\left[s_{1-n}(x',\xi')\right]\right\}\sigma(\xi')dx',
\end{align}
Then~~ a) ${\rm \widetilde{Wres}}([A,B])=0 $, for any
$A,B\in\mathcal{B}$;~~ b) It is a unique continuous trace on
$\mathcal{B}/\mathcal{B}^{-\infty}$.
\end{thm}

\begin{defn}{\rm\cite{Wa3} }
Lower dimensional volumes of spin manifolds with boundary are defined by
 \begin{equation}
{\rm Vol}^{(p_1,p_2)}_nM:= \widetilde{{\rm Wres}}[\pi^+D^{-p_1}\circ\pi^+D^{-p_2}],
\end{equation}
\end{defn}
 By \cite{Wa3}, we get
\begin{equation}
\widetilde{{\rm Wres}}[\pi^+D^{-p_1}\circ\pi^+D^{-p_2}]=\int_M\int_{|\xi|=1}{\rm
trace}_{\wedge^*T^*M\bigotimes\mathbb{C}}[\sigma_{-n}(D^{-p_1-p_2})]\sigma(\xi)dx+\int_{\partial M}\Phi,
\end{equation}
and
\begin{eqnarray}
\Phi&=\int_{|\xi'|=1}\int^{+\infty}_{-\infty}\sum^{\infty}_{j, k=0}\sum\frac{(-i)^{|\alpha|+j+k+1}}{\alpha!(j+k+1)!}
\times {\rm trace}_{\wedge^*T^*M\bigotimes\mathbb{C}}[\partial^j_{x_n}\partial^\alpha_{\xi'}\partial^k_{\xi_n}\sigma^+_{r}(D^{-p_1})(x',0,\xi',\xi_n)
\nonumber\\
&\times\partial^\alpha_{x'}\partial^{j+1}_{\xi_n}\partial^k_{x_n}\sigma_{l}(D^{-p_2})(x',0,\xi',\xi_n)]d\xi_n\sigma(\xi')dx',
\end{eqnarray}
 where the sum is taken over $r+l-k-|\alpha|-j-1=-n,~~r\leq -p_1,l\leq -p_2$.

 Since $[\sigma_{-n}(D^{-p_1-p_2})]|_M$ has the same expression as $\sigma_{-n}(D^{-p_1-p_2})$ in the case of manifolds without
boundary, so locally we can compute the first term by \cite{Ka}, \cite{KW}, \cite{Wa3}, \cite{Po}.

For any fixed point $x_0\in\partial M$, we choose the normal coordinates
$U$ of $x_0$ in $\partial M$ (not in $M$) and compute $\Phi(x_0)$ in the coordinates $\widetilde{U}=U\times [0,1)\subset M$ and the
metric $\frac{1}{h(x_n)}g^{\partial M}+dx_n^2.$ The dual metric of $g^M$ on $\widetilde{U}$ is ${h(x_n)}g^{\partial M}+dx_n^2.$  Write
$g^M_{ij}=g^M(\frac{\partial}{\partial x_i},\frac{\partial}{\partial x_j});~ g_M^{ij}=g^M(dx_i,dx_j)$, then

\begin{equation}
[g^M_{i,j}]= \left[\begin{array}{lcr}
  \frac{1}{h(x_n)}[g_{i,j}^{\partial M}]  & 0  \\
   0  &  1
\end{array}\right];~~~
[g_M^{i,j}]= \left[\begin{array}{lcr}
  h(x_n)[g^{i,j}_{\partial M}]  & 0  \\
   0  &  1
\end{array}\right],
\end{equation}
and
\begin{equation}
\partial_{x_s}g_{ij}^{\partial M}(x_0)=0, 1\leq i,j\leq n-1; ~~~g_{ij}^M(x_0)=\delta_{ij}.
\end{equation}
\indent From \cite{Wa3}, we can get three lemmas.
\begin{lem}{\rm \cite{Wa3}}\label{le:32}
With the metric $g^{M}$ on $M$ near the boundary
\begin{eqnarray}
\partial_{x_j}(|\xi|_{g^M}^2)(x_0)&=&\left\{
       \begin{array}{c}
        0,  ~~~~~~~~~~ ~~~~~~~~~~ ~~~~~~~~~~~~~{\rm if }~j<n, \\[2pt]
       h'(0)|\xi'|^{2}_{g^{\partial M}},~~~~~~~~~~~~~~~~~~~~{\rm if }~j=n;
       \end{array}
    \right. \\
\partial_{x_j}[c(\xi)](x_0)&=&\left\{
       \begin{array}{c}
      0,  ~~~~~~~~~~ ~~~~~~~~~~ ~~~~~~~~~~~~~{\rm if }~j<n,\\[2pt]
\partial x_{n}(c(\xi'))(x_{0}), ~~~~~~~~~~~~~~~~~{\rm if }~j=n,
       \end{array}
    \right.
\end{eqnarray}
where $\xi=\xi'+\xi_{n}dx_{n}$.
\end{lem}
\begin{lem}{\rm \cite{Wa3}}\label{le:32}With the metric $g^{M}$ on $M$ near the boundary
\begin{align}
\omega_{s,t}(e_i)(x_0)&=\left\{
       \begin{array}{c}
        \omega_{n,i}(e_i)(x_0)=\frac{1}{2}h'(0),  ~~~~~~~~~~ ~~~~~~~~~~~{\rm if }~s=n,t=i,i<n, \\[2pt]
       \omega_{i,n}(e_i)(x_0)=-\frac{1}{2}h'(0),~~~~~~~~~~~~~~~~~~~{\rm if }~s=i,t=n,i<n,\\[2pt]
    \omega_{s,t}(e_i)(x_0)=0,~~~~~~~~~~~~~~~~~~~~~~~~~~~other~cases,~~~~~~~~~\\[2pt]
       \end{array}
    \right.
\end{align}
where $(\omega_{s,t})$ denotes the connection matrix of Levi-Civita connection $\nabla^L$.
\end{lem}
\begin{lem}{\rm \cite{Wa3}}
\begin{align}
\Gamma_{st}^k(x_0)&=\left\{
       \begin{array}{c}
        \Gamma^n_{ii}(x_0)=\frac{1}{2}h'(0)  ~~~~~~~~~~ ~~~~~~~~~~~{\rm if }~s=t=i,k=n,i<n, \\[2pt]
        \Gamma^i_{ni}(x_0)=-\frac{1}{2}h'(0),~~~~~~~~~~~~~~~~~~~{\rm if }~s=n,t=i,k=i,i<n,\\[2pt]
        \Gamma^i_{in}(x_0)=-\frac{1}{2}h'(0),~~~~~~~~~~~~~~~~~~~{\rm if }~s=i,t=n,k=i,i<n,\\[2pt]
        \Gamma_{st}^i(x_0)=0,~~~~~~~~~~~~~~~~~~~~~~~~~~~other~cases.~~~~~~~~~
       \end{array}
    \right.
\end{align}
\end{lem}
\indent By (3.6) and (3.7), we firstly compute
\begin{equation}
\widetilde{{\rm Wres}}[\pi^+{\widetilde{D}}^{-1}\circ\pi^+({\widetilde{D}}^*)^{-1}]=\int_M\int_{|\xi|=1}{\rm
trace}_{\wedge^*T^*M\bigotimes\mathbb{C}}[\sigma_{-4}(({\widetilde{D}}^*{\widetilde{D}})^{-1})]\sigma(\xi)dx+\int_{\partial M}\Phi,
\end{equation}
where
\begin{align}
\Phi &=\int_{|\xi'|=1}\int^{+\infty}_{-\infty}\sum^{\infty}_{j, k=0}\sum\frac{(-i)^{|\alpha|+j+k+1}}{\alpha!(j+k+1)!}
\times {\rm trace}_{\wedge^*T^*M\bigotimes\mathbb{C}}[\partial^j_{x_n}\partial^\alpha_{\xi'}\partial^k_{\xi_n}\sigma^+_{r}({\widetilde{D}}^{-1})\nonumber\\
&(x',0,\xi',\xi_n)\times\partial^\alpha_{x'}\partial^{j+1}_{\xi_n}\partial^k_{x_n}\sigma_{l}(({\widetilde{D}}^*)^{-1})(x',0,\xi',\xi_n)]d\xi_n\sigma(\xi')dx',\nonumber\\
\end{align}
and the sum is taken over $r+l-k-j-|\alpha|=-3,~~r\leq -1,l\leq-1$.\\

\indent By Theorem 2.2, we can compute the interior of $\widetilde{{\rm Wres}}[\pi^+{\widetilde{D}}^{-1}\circ\pi^+({\widetilde{D}}^*)^{-1}]$, so
\begin{align}
&\int_M\int_{|\xi|=1}{\rm
trace}_{\wedge^*T^*M}[\sigma_{-4}(({\widetilde{D}}^*\widetilde{D})^{-1})]\sigma(\xi)dx\nonumber\\
&=32\pi^2\int_{M}\bigg\{-\frac{1}{3}s-12f_2\overline{f_2}+4f_2^2+4\overline{f_2}^2+4f_1\overline{f_1}\sum_{u<v}(p_{uv}-p_{vu})^2\bigg\}d{\rm Vol_{M}}.\nonumber\\
\end{align}

\indent Now we  need to compute $\int_{\partial M} \Phi$. Since, some operators have the following symbols.
\begin{lem} The following identities hold:
\begin{align}
\sigma_1({\widetilde{D}})&=\sigma_1({\widetilde{D}}^*)=ic(\xi); \nonumber\\ \sigma_0({\widetilde{D}})&=
-\frac{1}{4}\sum_{i,s,t}\omega_{s,t}(e_i)c(e_i)c(e_s)c(e_t)+(f_1\sum_{u<v}(p_{uv}-p_{vu})c(e_u)c(e_v)+f_2); \nonumber\\
\sigma_0({\widetilde{D}}^*)&=
-\frac{1}{4}\sum_{i,s,t}\omega_{s,t}(e_i)c(e_i)c(e_s)c(e_t)+(-\overline{f_1}\sum_{u<v}(p_{uv}-p_{vu})c(e_u)c(e_v)+\overline{f_2}).\nonumber\\
\end{align}
\end{lem}

\indent Write
 \begin{eqnarray}
D_x^{\alpha}&=(-i)^{|\alpha|}\partial_x^{\alpha};
~\sigma(\widetilde{D})=p_1+p_0;
~(\sigma(\widetilde{D})^{-1})=\sum^{\infty}_{j=1}q_{-j}.
\end{eqnarray}

\indent By the composition formula of pseudodifferential operators, we have
\begin{align}
1=\sigma(\widetilde{D}\circ {\widetilde{D}}^{-1})&=\sum_{\alpha}\frac{1}{\alpha!}\partial^{\alpha}_{\xi}[\sigma({\widetilde{D}})]
{\widetilde{D}}_x^{\alpha}[\sigma({\widetilde{D}}^{-1})]\nonumber\\
&=(p_1+p_0)(q_{-1}+q_{-2}+q_{-3}+\cdots)\nonumber\\
&~~~+\sum_j(\partial_{\xi_j}p_1+\partial_{\xi_j}p_0)(
D_{x_j}q_{-1}+D_{x_j}q_{-2}+D_{x_j}q_{-3}+\cdots)\nonumber\\
&=p_1q_{-1}+(p_1q_{-2}+p_0q_{-1}+\sum_j\partial_{\xi_j}p_1D_{x_j}q_{-1})+\cdots,
\end{align}
so
\begin{equation}
q_{-1}=p_1^{-1};~q_{-2}=-p_1^{-1}[p_0p_1^{-1}+\sum_j\partial_{\xi_j}p_1D_{x_j}(p_1^{-1})].
\end{equation}
\begin{lem} The following identities hold:
\begin{align}
\sigma_{-1}({\widetilde{D}}^{-1})&=\sigma_{-1}(({\widetilde{D}}^*)^{-1})=\frac{ic(\xi)}{|\xi|^2};\nonumber\\
\sigma_{-2}({\widetilde{D}}^{-1})&=\frac{c(\xi)\sigma_{0}({\widetilde{D}}^{-1})c(\xi)}{|\xi|^4}+\frac{c(\xi)}{|\xi|^6}\sum_jc(dx_j)
\Big[\partial_{x_j}(c(\xi))|\xi|^2-c(\xi)\partial_{x_j}(|\xi|^2)\Big] ;\nonumber\\
\sigma_{-2}(({\widetilde{D}}^*)^{-1})&=\frac{c(\xi)\sigma_{0}(({\widetilde{D}}^*)^{-1})c(\xi)}{|\xi|^4}+\frac{c(\xi)}{|\xi|^6}\sum_jc(dx_j)
\Big[\partial_{x_j}(c(\xi))|\xi|^2-c(\xi)\partial_{x_j}(|\xi|^2)\Big].\nonumber\\
\end{align}
\end{lem}
\indent When $n=4$, then ${\rm tr}_{\wedge^*T^*M}[{\rm \texttt{id}}]={\rm dim}(\wedge^*(4))=4$, where ${\rm tr}$ as shorthand of ${\rm trace}$, the sum is taken over $
r+l-k-j-|\alpha|=-3,~~r\leq -1,l\leq-1,$ then we have the following five cases:
~\\
\noindent  {\bf case a)~I)}~$r=-1,~l=-1,~k=j=0,~|\alpha|=1$\\

\noindent By (3.18), we get
\begin{equation}
\Phi_1=-\int_{|\xi'|=1}\int^{+\infty}_{-\infty}\sum_{|\alpha|=1}
 {\rm tr}[\partial^\alpha_{\xi'}\pi^+_{\xi_n}\sigma_{-1}({\widetilde{D}}^{-1})\times
 \partial^\alpha_{x'}\partial_{\xi_n}\sigma_{-1}(({\widetilde{D}}^*)^{-1})](x_0)d\xi_n\sigma(\xi')dx'.
\end{equation}
By Lemma 3.3, for $i<n$, then
\begin{equation}\partial_{x_i}\left(\frac{ic(\xi)}{|\xi|^2}\right)(x_0)=
\frac{i\partial_{x_i}[c(\xi)](x_0)}{|\xi|^2}
-\frac{ic(\xi)\partial_{x_i}(|\xi|^2)(x_0)}{|\xi|^4}=0,
\end{equation}
\noindent so $\Phi_1=0$.\\

 \noindent  {\bf case a)~II)}~$r=-1,~l=-1,~k=|\alpha|=0,~j=1$\\

\noindent By (3.18), we get
\begin{equation}
\Phi_2=-\frac{1}{2}\int_{|\xi'|=1}\int^{+\infty}_{-\infty} {\rm
trace} [\partial_{x_n}\pi^+_{\xi_n}\sigma_{-1}({\widetilde{D}}^{-1})\times
\partial_{\xi_n}^2\sigma_{-1}(({\widetilde{D}}^*)^{-1})](x_0)d\xi_n\sigma(\xi')dx'.
\end{equation}
\noindent By Lemma 3.7, we have\\
\begin{eqnarray}\partial^2_{\xi_n}\sigma_{-1}(({\widetilde{D}}^*)^{-1})(x_0)=i\left(-\frac{6\xi_nc(dx_n)+2c(\xi')}
{|\xi|^4}+\frac{8\xi_n^2c(\xi)}{|\xi|^6}\right);
\end{eqnarray}
\begin{eqnarray}
\partial_{x_n}\sigma_{-1}({\widetilde{D}}^{-1})(x_0)=\frac{i\partial_{x_n}c(\xi')(x_0)}{|\xi|^2}-\frac{ic(\xi)|\xi'|^2h'(0)}{|\xi|^4}.
\end{eqnarray}
By (3.5), (3.6), we get
\begin{align}
\pi^+_{\xi_n}\left[\frac{c(\xi)}{|\xi|^4}\right](x_0)|_{|\xi'|=1}&=\pi^+_{\xi_n}\left[\frac{c(\xi')+\xi_nc(dx_n)}{(1+\xi_n^2)^2}\right]\nonumber\\
&=\frac{1}{2\pi i}{\rm lim}_{u\rightarrow
0^-}\int_{\Gamma^+}\frac{\frac{c(\xi')+\eta_nc(dx_n)}{(\eta_n+i)^2(\xi_n+iu-\eta_n)}}
{(\eta_n-i)^2}d\eta_n\nonumber\\
&=-\frac{(i\xi_n+2)c(\xi')+ic(dx_n)}{4(\xi_n-i)^2}.
\end{align}
Similarly we have,
\begin{eqnarray}
\pi^+_{\xi_n}\left[\frac{i\partial_{x_n}c(\xi')}{|\xi|^2}\right](x_0)|_{|\xi'|=1}=\frac{\partial_{x_n}[c(\xi')](x_0)}{2(\xi_n-i)}.
\end{eqnarray}
By (3.29), then\\
\begin{align}\pi^+_{\xi_n}\partial_{x_n}\sigma_{-1}({\widetilde{D}}^{-1})|_{|\xi'|=1}
=\frac{\partial_{x_n}[c(\xi')](x_0)}{2(\xi_n-i)}+ih'(0)
\left[\frac{(i\xi_n+2)c(\xi')+ic(dx_n)}{4(\xi_n-i)^2}\right].
\end{align}
\noindent By the relation of the Clifford action and ${\rm tr}{AB}={\rm tr }{BA}$, we have the equalities:\\
\begin{align}
&{\rm tr}[c(\xi')c(dx_n)]=0;~~{\rm tr}[c(dx_n)^2]=-4;~~{\rm tr}[c(\xi')^2](x_0)|_{|\xi'|=1}=-4;\nonumber\\
&{\rm tr}[\partial_{x_n}c(\xi')c(dx_n)]=0;~~{\rm tr}[\partial_{x_n}c(\xi')c(\xi')](x_0)|_{|\xi'|=1}=-2h'(0).\nonumber\\
\end{align}
By (3.31), we have
\begin{eqnarray}
&h'(0){\rm tr}\bigg[\frac{(i\xi_n+2)c(\xi')+ic(dx_n)}{4(\xi_n-i)^2}\times
\bigg(\frac{6\xi_nc(dx_n)+2c(\xi')}{(1+\xi_n^2)^2}
-\frac{8\xi_n^2[c(\xi')+\xi_nc(dx_n)]}{(1+\xi_n^2)^3}\bigg)
\bigg](x_0)|_{|\xi'|=1}\nonumber\\
&=-4h'(0)\frac{-2i\xi_n^2-\xi_n+i}{(\xi_n-i)^4(\xi_n+i)^3}.
\end{eqnarray}
Similarly, we have
\begin{align}
&-i{\rm
tr}\bigg[\bigg(\frac{\partial_{x_n}[c(\xi')](x_0)}{2(\xi_n-i)}\bigg)
\times\bigg(\frac{6\xi_nc(dx_n)+2c(\xi')}{(1+\xi_n^2)^2}-\frac{8\xi_n^2[c(\xi')+\xi_nc(dx_n)]}
{(1+\xi_n^2)^3}\bigg)\bigg](x_0)|_{|\xi'|=1}\nonumber\\
&=-2ih'(0)\frac{3\xi_n^2-1}{(\xi_n-i)^4(\xi_n+i)^3}.\nonumber\\
\end{align}
Then\\
\begin{align}
\Phi_2&=-\int_{|\xi'|=1}\int^{+\infty}_{-\infty}\frac{ih'(0)(\xi_n-i)^2}
{(\xi_n-i)^4(\xi_n+i)^3}d\xi_n\sigma(\xi')dx'\\\notag
&=-ih'(0)\Omega_3\int_{\Gamma^+}\frac{1}{(\xi_n-i)^2(\xi_n+i)^3}d\xi_ndx'\\\notag
&=-ih'(0)\Omega_32\pi i[\frac{1}{(\xi_n+i)^3}]^{(1)}|_{\xi_n=i}dx'\\\notag
&=-\frac{3}{8}\pi h'(0)\Omega_3dx',
\end{align}
where ${\rm \Omega_{3}}$ is the canonical volume of $S^{3}.$\\

\noindent  {\bf case a)~III)}~$r=-1,~l=-1,~j=|\alpha|=0,~k=1$\\

\noindent By (3.18), we get
\begin{equation}
\Phi_3=-\frac{1}{2}\int_{|\xi'|=1}\int^{+\infty}_{-\infty}
{\rm trace} [\partial_{\xi_n}\pi^+_{\xi_n}\sigma_{-1}({\widetilde{D}}^{-1})\times
\partial_{\xi_n}\partial_{x_n}\sigma_{-1}(({\widetilde{D}}^*)^{-1})](x_0)d\xi_n\sigma(\xi')dx'.
\end{equation}
\noindent By Lemma 3.7, we have\\
\begin{equation}
\partial_{\xi_n}\partial_{x_n}\sigma_{-1}(({\widetilde{D}}^*)^{-1})(x_0)|_{|\xi'|=1}
=-ih'(0)\left[\frac{c(dx_n)}{|\xi|^4}-4\xi_n\frac{c(\xi')
+\xi_nc(dx_n)}{|\xi|^6}\right]-\frac{2\xi_ni\partial_{x_n}c(\xi')(x_0)}{|\xi|^4};
\end{equation}
\begin{equation}
\partial_{\xi_n}\pi^+_{\xi_n}\sigma_{-1}({\widetilde{D}}^{-1})(x_0)|_{|\xi'|=1}=-\frac{c(\xi')+ic(dx_n)}{2(\xi_n-i)^2}.
\end{equation}
Similar to {\rm case~a)~II)}, we have\\
\begin{equation}
{\rm tr}\left\{\frac{c(\xi')+ic(dx_n)}{2(\xi_n-i)^2}\times
ih'(0)\left[\frac{c(dx_n)}{|\xi|^4}-4\xi_n\frac{c(\xi')+\xi_nc(dx_n)}{|\xi|^6}\right]\right\}=2h'(0)\frac{i-3\xi_n}{(\xi_n-i)^4(\xi_n+i)^3}
\end{equation}
and
\begin{eqnarray}
{\rm tr}\left[\frac{c(\xi')+ic(dx_n)}{2(\xi_n-i)^2}\times
\frac{2\xi_ni\partial_{x_n}c(\xi')(x_0)}{|\xi|^4}\right]
=\frac{-2ih'(0)\xi_n}{(\xi_n-i)^4(\xi_n+i)^2}.
\end{eqnarray}
So we have
\begin{align}
\Phi_3&=-\int_{|\xi'|=1}\int^{+\infty}_{-\infty}\frac{h'(0)(i-3\xi_n)}
{(\xi_n-i)^4(\xi_n+i)^3}d\xi_n\sigma(\xi')dx'\nonumber\\
&-\int_{|\xi'|=1}\int^{+\infty}_{-\infty}\frac{h'(0)i\xi_n}
{(\xi_n-i)^4(\xi_n+i)^2}d\xi_n\sigma(\xi')dx'\nonumber\\
&=-h'(0)\Omega_3\frac{2\pi i}{3!}[\frac{(i-3\xi_n)}{(\xi_n+i)^3}]^{(3)}|_{\xi_n=i}dx'+h'(0)\Omega_3\frac{2\pi i}{3!}[\frac{i\xi_n}{(\xi_n+i)^2}]^{(3)}|_{\xi_n=i}dx'\nonumber\\
&=\frac{3}{8}\pi h'(0)\Omega_3dx'.\nonumber\\
\end{align}

\noindent  {\bf case b)}~$r=-2,~l=-1,~k=j=|\alpha|=0$\\

\noindent By (3.18), we get
\begin{align}
\Phi_4&=-i\int_{|\xi'|=1}\int^{+\infty}_{-\infty}{\rm trace} [\pi^+_{\xi_n}\sigma_{-2}({\widetilde{D}}^{-1})\times
\partial_{\xi_n}\sigma_{-1}(({\widetilde{D}}^*)^{-1})](x_0)d\xi_n\sigma(\xi')dx'.
\end{align}
 By Lemma 3.7 we have\\
\begin{align}
\sigma_{-2}({\widetilde{D}}^{-1})(x_0)=\frac{c(\xi)\sigma_{0}({\widetilde{D}})(x_0)c(\xi)}{|\xi|^4}+\frac{c(\xi)}{|\xi|^6}c(dx_n)
[\partial_{x_n}[c(\xi')](x_0)|\xi|^2-c(\xi)h'(0)|\xi|^2_{\partial
M}],
\end{align}
where
\begin{align}
\sigma_{0}({\widetilde{D}})(x_0)&=-\frac{1}{4}\sum_{s,t,i}\omega_{s,t}(e_i)
(x_{0})c(e_i)c(e_s)c(e_t)\nonumber\\
&+f_1\sum_{u<v}(p_{uv}-p_{vu})c(e_u)c(e_v)+f_2.
\end{align}
We denote
\begin{align}
Q(x_0)&=-\frac{1}{4}\sum_{s,t,i}\omega_{s,t}(e_i)
(x_{0})c(e_i)c(e_s)c(e_t).
\end{align}
Then
\begin{align}
&\pi^+_{\xi_n}\sigma_{-2}({\widetilde{D}}^{-1}(x_0))|_{|\xi'|=1}\nonumber\\
&=\pi^+_{\xi_n}
\Big[\frac{c(\xi)(f_1\sum_{u<v}(p_{uv}-p_{vu})c(e_u)c(e_v)+f_2)(x_0)c(\xi)}{(1+\xi_n^2)^2}\Big]
\nonumber\\
&+\pi^+_{\xi_n}\Big[\frac{c(\xi)Q(x_0)c(\xi)+c(\xi)c(dx_n)\partial_{x_n}[c(\xi')](x_0)}{(1+\xi_n^2)^2}-h'(0)\frac{c(\xi)c(dx_n)c(\xi)}{(1+\xi_n^{2})^3}\Big].
\end{align}
By computations, we have the equalities:\\
\begin{align}
&{\rm tr}[(f_1\sum_{u<v}(p_{uv}-p_{vu})c(e_u)c(e_v)+f_2)c(dx_n)]=0;~~~\nonumber\\
&{\rm tr}[(f_1\sum_{u<v}(p_{uv}-p_{vu})c(e_u)c(e_v)+f_2)c(\xi')]=0.\nonumber\\
\end{align}
Since
\begin{align}
\partial_{\xi_n}\sigma_{-1}(({\widetilde{D}}^*)^{-1})=\partial_{\xi_n}q_{-1}(x_0)|_{|\xi'|=1}=i\left[\frac{c(dx_n)}{1+\xi_n^2}-\frac{2\xi_nc(\xi')+2\xi_n^2c(dx_n)}{(1+\xi_n^2)^2}\right].
\end{align}
Then, we have
\begin{equation}
\pi^+_{\xi_n}\Big[\frac{c(\xi)Q(x_0)c(\xi)+c(\xi)c(dx_n)\partial_{x_n}[c(\xi')](x_0)}{(1+\xi_n^2)^2}\Big]-h'(0)\pi^+_{\xi_n}\Big[\frac{c(\xi)c(dx_n)c(\xi)}{(1+\xi_n)^3}\Big]:= C_1-C_2,
\end{equation}
where
\begin{align}
C_1&=\frac{-1}{4(\xi_n-i)^2}[(2+i\xi_n)c(\xi')Q(x_0)c(\xi')+i\xi_nc(dx_n)Q(x_0)c(dx_n)\nonumber\\
&+(2+i\xi_n)c(\xi')c(dx_n)\partial_{x_n}c(\xi')+ic(dx_n)Q(x_0)c(\xi')
+ic(\xi')Q(x_0)c(dx_n)-i\partial_{x_n}c(\xi')]\nonumber\\
\end{align}
and
\begin{align}
C_2&=\frac{h'(0)}{2}\left[\frac{c(dx_n)}{4i(\xi_n-i)}+\frac{c(dx_n)-ic(\xi')}{8(\xi_n-i)^2}
+\frac{3\xi_n-7i}{8(\xi_n-i)^3}[ic(\xi')-c(dx_n)]\right].
\end{align}
By (3.50) and (3.52), we have\\
\begin{eqnarray}{\rm tr }[C_2\times\partial_{\xi_n}\sigma_{-1}(({\widetilde{D}}^*)^{-1})]|_{|\xi'|=1}
&=\frac{i}{2}h'(0)\frac{-i\xi_n^2-\xi_n+4i}{4(\xi_n-i)^3(\xi_n+i)^2}{\rm tr}[ \texttt{id}]\nonumber\\
&=2ih'(0)\frac{-i\xi_n^2-\xi_n+4i}{4(\xi_n-i)^3(\xi_n+i)^2}.
\end{eqnarray}
By (3.50) and (3.51), we have
\begin{eqnarray}{\rm tr }[C_1\times\partial_{\xi_n}\sigma_{-1}(({\widetilde{D}}^*)^{-1})]|_{|\xi'|=1}=
\frac{-2ic_0}{(1+\xi_n^2)^2}+h'(0)\frac{\xi_n^2-i\xi_n-2}{2(\xi_n-i)(1+\xi_n^2)^2},
\end{eqnarray}
where $Q=c_0c(dx_n)$ and $c_0=-\frac{3}{4}h'(0)$.\\

By (3.53) and (3.54), we have
\begin{align}
&-i\int_{|\xi'|=1}\int^{+\infty}_{-\infty}{\rm trace} [(C_1-C_2)\times
\partial_{\xi_n}\sigma_{-1}(({\widetilde{D}}^*)^{-1})](x_0)d\xi_n\sigma(\xi')dx'\nonumber\\
&=-\Omega_3\int_{\Gamma^+}\frac{2c_0(\xi_n-i)+ih'(0)}{(\xi_n-i)^3(\xi_n+i)^2}d\xi_ndx'\nonumber\\
&=\frac{9}{8}\pi h'(0)\Omega_3dx'.
\end{align}
Then, we have
\begin{align}
&-i\int_{|\xi'|=1}\int^{+\infty}_{-\infty}{\rm trace} [\pi^+_{\xi_n}
\Big[\frac{c(\xi)(f_1\sum_{u<v}(p_{uv}-p_{vu})c(e_u)c(e_v)+f_2)c(\xi)}{(1+\xi_n^2)^2}\Big]\nonumber\\
&\times
\partial_{\xi_n}\sigma_{-1}(({\widetilde{D}}^*)^{-1})](x_0)d\xi_n\sigma(\xi')dx'\nonumber\\
&=\frac{\pi}{4}{\rm tr}[c(dx_n)(f_1\sum_{u<v}(p_{uv}-p_{vu})c(e_u)c(e_v)+f_2)]\Omega_3dx'\nonumber\\
&=0.\nonumber\\
\end{align}
Then, we have\\
\begin{align}
\Phi_4=\frac{9}{8}\pi h'(0)\Omega_3dx'.
\end{align}

\noindent {\bf  case c)}~$r=-1,~l=-2,~k=j=|\alpha|=0$\\
By (3.18), we get
\begin{align}
\Phi_5=-i\int_{|\xi'|=1}\int^{+\infty}_{-\infty}{\rm trace} [\pi^+_{\xi_n}\sigma_{-1}({\widetilde{D}}^{-1})\times
\partial_{\xi_n}\sigma_{-2}(({\widetilde{D}}^*)^{-1})](x_0)d\xi_n\sigma(\xi')dx'.
\end{align}
By (3.5) and (3.6), Lemma 3.7, we have
\begin{align}
\pi^+_{\xi_n}\sigma_{-1}({\widetilde{D}}^{-1})|_{|\xi'|=1}=\frac{c(\xi')+ic(dx_n)}{2(\xi_n-i)}.
\end{align}

Since
\begin{equation}
\sigma_{-2}(({\widetilde{D}}^*)^{-1})(x_0)=\frac{c(\xi)\sigma_{0}({\widetilde{D}}^*)(x_0)c(\xi)}{|\xi|^4}+\frac{c(\xi)}{|\xi|^6}c(dx_n)
\bigg[\partial_{x_n}[c(\xi')](x_0)|\xi|^2-c(\xi)h'(0)|\xi|^2_{\partial_
M}\bigg],
\end{equation}
where
\begin{align}
\sigma_{0}({\widetilde{D}}^*)(x_0)&=
-\frac{1}{4}\sum_{s,t,i}\omega_{s,t}(e_i)(x_{0})c(e_i)c(e_s)c(e_t)\nonumber\\
&+(-\overline{f_1}\sum_{u<v}(p_{uv}-p_{vu})c(e_u)c(e_v)+\overline{f_2})(x_{0})\nonumber\\
&=Q(x_0)+(-\overline{f_1}\sum_{u<v}(p_{uv}-p_{vu})c(e_u)c(e_v)+\overline{f_2})(x_{0}),
\end{align}
then
\begin{align}
\partial_{\xi_n}\sigma_{-2}(({\widetilde{D}}^*)^{-1})(x_0)|_{|\xi'|=1}&=
\partial_{\xi_n}\bigg\{\frac{c(\xi)[Q(x_0)
+(-\overline{f_1}\sum_{u<v}(p_{uv}-p_{vu})c(e_u)c(e_v)+\overline{f_2})(x_{0})]c(\xi)}{|\xi|^4}\nonumber\\
&+\frac{c(\xi)}{|\xi|^6}c(dx_n)[\partial_{x_n}[c(\xi')](x_0)|\xi|^2-c(\xi)h'(0)]\bigg\}\nonumber\\
&=\partial_{\xi_n}\bigg\{\frac{c(\xi)}{|\xi|^6}c(dx_n)[\partial_{x_n}[c(\xi')](x_0)|\xi|^2-c(\xi)h'(0)]\bigg\}+\partial_{\xi_n}\frac{c(\xi)Q(x_0)c(\xi)}{|\xi|^4}\nonumber\\
&+\partial_{\xi_n}\frac{c(\xi)(-\overline{f_1}\sum_{u<v}(p_{uv}-p_{vu})c(e_u)c(e_v)+\overline{f_2})(x_{0})c(\xi)}{|\xi|^4}.\nonumber\\
\end{align}
By computations, we have
\begin{align}
&\partial_{\xi_n}\frac{c(\xi)(-\overline{f_1}\sum_{u<v}(p_{uv}-p_{vu})c(e_u)c(e_v)+\overline{f_2})(x_{0})c(\xi)}{|\xi|^4}\nonumber\\
&=\frac{c(dx_n)(-\overline{f_1}\sum_{u<v}(p_{uv}-p_{vu})c(e_u)c(e_v)+\overline{f_2})(x_{0})c(\xi)}{|\xi|^4}\nonumber\\
&+\frac{c(\xi)(-\overline{f_1}\sum_{u<v}(p_{uv}-p_{vu})c(e_u)c(e_v)+\overline{f_2})(x_{0})c(dx_n)}{|\xi|^4}\nonumber\\
&-\frac{4\xi_n c(\xi)(-\overline{f_1}\sum_{u<v}(p_{uv}-p_{vu})c(e_u)c(e_v)+\overline{f_2})(x_{0})c(\xi)}{|\xi|^4}.\nonumber\\
\end{align}
We denote $$q_{-2}^{1}=\frac{c(\xi)Q(x_0)c(\xi)}{|\xi|^4}+\frac{c(\xi)}{|\xi|^6}c(dx_n)[\partial_{x_n}[c(\xi')](x_0)|\xi|^2-c(\xi)h'(0)],$$ then
\begin{align}
\partial_{\xi_n}(q_{-2}^{1})&=\frac{1}{(1+\xi_n^2)^3}\bigg[(2\xi_n-2\xi_n^3)c(dx_n)Q(x_0)c(dx_n)
+(1-3\xi_n^2)c(dx_n)Q(x_0)c(\xi')\nonumber\\
&+(1-3\xi_n^2)c(\xi')Q(x_0)c(dx_n)
-4\xi_nc(\xi')Q(x_0)c(\xi')
+(3\xi_n^2-1)\partial_{x_n}c(\xi')\nonumber\\
&-4\xi_nc(\xi')c(dx_n)\partial_{x_n}c(\xi')
+2h'(0)c(\xi')+2h'(0)\xi_nc(dx_n)\bigg]\nonumber\\
&+6\xi_nh'(0)\frac{c(\xi)c(dx_n)c(\xi)}{(1+\xi^2_n)^4}.
\end{align}
By (3.59) and (3.64), we have
\begin{equation}
{\rm tr}[\pi^+_{\xi_n}\sigma_{-1}({\widetilde{D}}^{-1})\times
\partial_{\xi_n}(q^1_{-2})](x_0)|_{|\xi'|=1}
=\frac{3h'(0)(i\xi^2_n+\xi_n-2i)}{(\xi-i)^3(\xi+i)^3}
+\frac{12h'(0)i\xi_n}{(\xi-i)^3(\xi+i)^4},
\end{equation}
then
\begin{equation}
-i\Omega_3\int_{\Gamma_+}[\frac{3h'(0)(i\xi_n^2+\xi_n-2i)}
{(\xi_n-i)^3(\xi_n+i)^3}+\frac{12h'(0)i\xi_n}{(\xi_n-i)^3(\xi_n+i)^4}]d\xi_ndx'=
-\frac{9}{8}\pi h'(0)\Omega_3dx'.
\end{equation}
By $\int_{|\xi'|=1}\{\xi_{1}\cdot\cdot\cdot\xi_{2d+1}\}\sigma(\xi')=0$, (3.59) and (3.62), we have
\begin{align}
&-i\int_{|\xi'|=1}\int^{+\infty}_{-\infty}{\rm tr}[\pi^+_{\xi_n}\sigma_{-1}({\widetilde{D}}^{-1})\times
\partial_{\xi_n}\frac{c(\xi)(-\overline{f_1}\sum_{u<v}(p_{uv}-p_{vu})c(e_u)c(e_v)+\overline{f_2})c(\xi)}
{|\xi|^4}]\nonumber\\
&(x_0)d\xi_n\sigma(\xi')dx'\nonumber\\
&=-i\int_{|\xi'|=1}\int^{+\infty}_{-\infty}\frac{i}{(\xi-i)(\xi+i)^3}{\rm tr}[c(dx_n)(-\overline{f_1}\sum_{u<v}(p_{uv}-p_{vu})c(e_u)c(e_v)+\overline{f_2})]\nonumber\\
&(x_0)d\xi_n\sigma(\xi')dx'\nonumber\\
&=-\frac{\pi}{4}{\rm tr}[c(dx_n)(-\overline{f_1}\sum_{u<v}(p_{uv}-p_{vu})c(e_u)c(e_v)+\overline{f_2})]\Omega_3dx'\nonumber\\
&=0.\nonumber\\
\end{align}
Then,
\begin{align}
\Phi_5=-\frac{9}{8}\pi h'(0)\Omega_3dx'.
\end{align}
So $\Phi=\Phi_1+\Phi_2+\Phi_3+\Phi_4+\Phi_5=0$.

\begin{thm}
Let $M$ be a $4$-dimensional
compact oriented spin manifolds with the boundary $\partial M$ and the metric
$g^M$ as above, ${\widetilde{D}}$ and ${\widetilde{D}}^*$ be the Dirac-Witten operators on $\widetilde{M}$, then
\begin{align}
&\widetilde{{\rm Wres}}[\pi^+{\widetilde{D}}^{-1}\circ\pi^+({\widetilde{D}}^*)^{-1}]\nonumber\\
&=32\pi^2\int_{M}\bigg\{-\frac{1}{3}s-12f_2\overline{f_2}+4f_2^2+4\overline{f_2}^2+4f_1\overline{f_1}\sum_{u<v}(p_{uv}-p_{vu})^2\bigg\}d{\rm Vol_{M}}.\nonumber\\
\end{align}
where $s$ is the scalar curvature.
\end{thm}

\section{A Kastler-Kalau-Walze type theorem for $6$-dimensional manifolds with boundary }
Firstly, we prove the Kastler-Kalau-Walze type theorems for $6$-dimensional manifolds with boundary. From \cite{Wa5}, we know that
\begin{equation}
\widetilde{{\rm Wres}}[\pi^+{\widetilde{D}}^{-1}\circ\pi^+({\widetilde{D}}^{*}{\widetilde{D}}
      {\widetilde{D}}^{*})^{-1}]=\int_M\int_{|\xi|=1}{\rm
trace}_{\wedge ^*T^*M\bigotimes\mathbb{C}}[\sigma_{-4}(({\widetilde{D}}^*{\widetilde{D}})^{-2})]\sigma(\xi)dx+\int_{{\partial_t}M}\Psi,
\end{equation}
where
\begin{align}
\Psi &=\int_{|\xi'|=1}\int^{+\infty}_{-\infty}\sum^{\infty}_{j, k=0}\sum\frac{(-i)^{|\alpha|+j+k+1}}{\alpha!(j+k+1)!}
\times {\rm trace}_{\wedge ^*T^*M\bigotimes\mathbb{C}}[\partial^j_{x_n}\partial^\alpha_{\xi'}{\partial_t}^k_{\xi_n}\sigma^+_{r}({\widetilde{D}}^{-1})(x',0,\xi',\xi_n)
\nonumber\\
&\times\partial^\alpha_{x'}\partial^{j+1}_{\xi_n}\partial^k_{x_n}\sigma_{l}
(({\widetilde{D}}^{*}{\widetilde{D}}
      {\widetilde{D}}^{*})^{-1})(x',0,\xi',\xi_n)]d\xi_n\sigma(\xi')dx',\nonumber\\
\end{align}
and the sum is taken over $r+\ell-k-j-|\alpha|-1=-6, \ r\leq-1, \ell\leq -3$.\\
By Theorem 2.2, we compute the interior term of (4.1), then
\begin{align}
&\int_M\int_{|\xi|=1}{\rm
trace}_{\wedge^*T^*M\bigotimes\mathbb{C}}[\sigma_{-4}(({\widetilde{D}}^*{\widetilde{D}})^{-2})]\sigma(\xi)dx\nonumber\\
&=128\pi^3\int_{M}\bigg(
-\frac{2}{3}s-4(f_1^2+\overline{f_1}^2-4f_1\overline{f_1})\sum_{u<v}(p_{uv}-p_{vu})^2-32f_2\overline{f_2}+12f_2^2+12\overline{f_2}^2\bigg)d{\rm Vol_{M}}.\nonumber\\
\end{align}

Next, we compute $\int_{\partial M} \Psi$. By computations, we get
\begin{align}
{\widetilde{D}}^*{\widetilde{D}}{\widetilde{D}}^*
&=\sum^{n}_{i=1}c(e_{i})\langle e_{i},dx_{l}\rangle(-g^{ij}\partial_{l}\partial_{i}\partial_{j})
+\sum^{n}_{i=1}c(e_{i})\langle e_{i},dx_{l}\rangle \bigg\{-(\partial_{l}g^{ij})\partial_{i}\partial_{j}-g^{ij}\bigg(4\sigma_{i}
\partial_{j}-2\Gamma^{k}_{ij}\partial_{k}\bigg)\partial_{l}\bigg\}\nonumber\\
&+\sum^{n}_{i=1}c(e_{i})\langle e_{i},dx_{l}\rangle \bigg\{-2(\partial_{l}g^{ij})\sigma_{i}\partial_{j}+g^{ij}
(\partial_{l}\Gamma^{k}_{ij})\partial_{k}-2g^{ij}\partial_{l}\sigma_{i}\partial_{j}
+(\partial_{l}g^{ij})\Gamma^{k}_{ij}\partial_{k}+\sum_{j,k}\Big[\partial_{l}\Big((f_1\nonumber\\
&\sum_{u<v}(p_{uv}-p_{vu})c(e_u)c(e_v)+f_2)c(e_{j})+c(e_{j})(\overline{f_1}\sum_{u<v}(p_{uv}-p_{vu})c(e_u)c(e_v)-\overline{f_2})\Big)\Big]
\langle      e_{j},dx^{k}\rangle\partial_{k}\nonumber\\
&+\sum_{j,k}\Big(f_1\sum_{u<v}(p_{uv}-p_{vu})c(e_u)c(e_v)+f_2)c(e_{j})+c(e_{j})(\overline{f_1}\sum_{u<v}(p_{uv}-p_{vu})c(e_u)c(e_v)-\overline{f_2})\Big)\nonumber\\
&\Big[\partial_{l}\langle e_{j},dx^{k}\rangle\Big]\partial_{k} \bigg\}+\sum^{n}_{i=1}c(e_{i})\langle e_{i},dx_{l}\rangle\partial_{l}\bigg\{-g^{ij}\Big[(\partial_{i}\sigma_{j})+\sigma_{i}\sigma_{j}-\Gamma_{i,j}^{k}\sigma_{k}+\sum_{i,j}g^{i,j}\Big[(f_1\sum_{u<v}(p_{uv}-p_{vu})\nonumber\\
&c(e_u)c(e_v)+f_2)c(\partial_{i})\sigma_{i}+c(\partial_{i})\partial_{i}(-\overline{f_1}\sum_{u<v}(p_{uv}-p_{vu})c(e_u)c(e_v)+\overline{f_2})+c(\partial_{i})\sigma_{i}(-\overline{f_1}\sum_{u<v}(p_{uv}-p_{vu})\nonumber\\
&c(e_u)c(e_v)+\overline{f_2})\Big]+\frac{1}{4}s+(\overline{f_1}\sum_{u<v}(p_{uv}-p_{vu})c(e_u)c(e_v)-\overline{f_2})^2\bigg\}+\Big[\sigma_{i}-(\overline{f_1}\sum_{u<v}(p_{uv}-p_{vu})c(e_u)\nonumber\\
&c(e_v)-\overline{f_2})\Big](-g^{ij}\partial_{i}\partial_{j})+\sum^{n}_{i=1}c(e_{i})\langle e_{i},dx_{l}\rangle \bigg\{2\sum_{j,k}\Big[(f_1\sum_{u<v}(p_{uv}-p_{vu})c(e_u)c(e_v)+f_2)c(e_{j})+c(e_{j})\nonumber\\
&(-\overline{f_1}\sum_{u<v}(p_{uv}-p_{vu})c(e_u)c(e_v)+\overline{f_2})\Big]\times\langle e_{i},dx_{k}\rangle\bigg\}_{l}\partial_{k}
+\Big[\sigma_{i}+(-\overline{f_1}\sum_{u<v}(p_{uv}-p_{vu})c(e_u)c(e_v)+\overline{f_2})\Big]\nonumber\\
&
\bigg\{-\sum_{i,j}g^{i,j}\Big[2\sigma_{i}\partial_{j}
-\Gamma_{i,j}^{k}\partial_{k}+(\partial_{i}\sigma_{j})+\sigma_{i}\sigma_{j} -\Gamma_{i,j}^{k}\sigma_{k}\Big]-\sum_{i,j}g^{i,j}\Big[c(\partial_{i})(-\overline{f_1}\sum_{u<v}(p_{uv}-p_{vu})c(e_u)c(e_v)\nonumber\\
&+\overline{f_2})+(f_1\sum_{u<v}(p_{uv}-p_{vu})c(e_u)c(e_v)+f_2)c(\partial_{i})\Big]\partial_{j}+\sum_{i,j}g^{i,j}\Big[(f_1\sum_{u<v}(p_{uv}-p_{vu})c(e_u)c(e_v)+f_2)c(\partial_{i})\nonumber\\
&\sigma_{i}+c(\partial_{i})\partial_{i}(-\overline{f_1}\sum_{u<v}(p_{uv}-p_{vu})c(e_u)c(e_v)+\overline{f_2})-c(\partial_{i})\sigma_{i}(-\overline{f_1}\sum_{u<v}(p_{uv}-p_{vu})c(e_u)c(e_v)+\overline{f_2})+c(\partial_{i})\nonumber\\
&\partial_{i}(-\overline{f_1}\sum_{u<v}(p_{uv}-p_{vu})c(e_u)c(e_v)+\overline{f_2})c(\partial_{i})\sigma_{i}(-\overline{f_1}\sum_{u<v}(p_{uv}-p_{vu})c(e_u)c(e_v)+\overline{f_2})\Big]+\frac{1}{4}s-(f_1\sum_{u<v}\nonumber\\
&(p_{uv}-p_{vu})c(e_u)c(e_v)+f_2)^2\bigg\}.\nonumber\\
\end{align}
Then, we obtain
\begin{lem} The following identities hold:
\begin{align}
\sigma_2({\widetilde{D}}^*{\widetilde{D}}{\widetilde{D}}^*)&=\sum_{i,j,l}c(dx_{l})\partial_{l}(g^{i,j})\xi_{i}\xi_{j} +c(\xi)(4\sigma^k-2\Gamma^k)\xi_{k}+2[|\xi|^2(-\overline{f_1}\sum_{u<v}(p_{uv}-p_{vu})\nonumber\\
&c(e_u)c(e_v)+\overline{f_2})-c(\xi)(f_1\sum_{u<v}(p_{uv}-p_{vu})c(e_u)c(e_v)+f_2)c(\xi)]\nonumber\\
&-\frac{1}{4}|\xi|^2\sum_{s,t,l}\omega_{s,t}
c(e_l)c(e_s)c(e_t)+|\xi|^2(-\overline{f_1}\sum_{u<v}(p_{uv}-p_{vu})c(e_u)c(e_v)+\overline{f_2})^2;\nonumber\\
\sigma_{3}({\widetilde{D}}^*{\widetilde{D}}{\widetilde{D}}^*)
&=ic(\xi)|\xi|^{2}.\nonumber\\
\end{align}
\end{lem}

Write
\begin{align}
\sigma({\widetilde{D}}^*{\widetilde{D}}{\widetilde{D}}^*)&=p_3+p_2+p_1+p_0;
~\sigma(({\widetilde{D}}^*{\widetilde{D}}{\widetilde{D}}^*)^{-1})=\sum^{\infty}_{j=3}q_{-j}.
\end{align}

By the composition formula of pseudodifferential operators, we have
\begin{align}
1=\sigma(({\widetilde{D}}^*{\widetilde{D}}{\widetilde{D}}^*)\circ ({\widetilde{D}}^*{\widetilde{D}}{\widetilde{D}}^*)^{-1})&=
\sum_{\alpha}\frac{1}{\alpha!}\partial^{\alpha}_{\xi}
[\sigma({\widetilde{D}}^*{\widetilde{D}}{\widetilde{D}}^*)]D^{\alpha}_{x}
[({\widetilde{D}}^*{\widetilde{D}}{\widetilde{D}}^*)^{-1}] \nonumber\\
&=(p_3+p_2+p_1+p_0)(q_{-3}+q_{-4}+q_{-5}+\cdots)\nonumber\\
&+\sum_j(\partial_{\xi_j}p_3+\partial_{\xi_j}p_2+\partial_{\xi_j}p_1+\partial_{\xi_j}p_0)\nonumber\\
&(D_{x_j}q_{-3}+D_{x_j}q_{-4}+D_{x_j}q_{-5}+\cdots) \nonumber\\
&=p_3q_{-3}+(p_3q_{-4}+p_2q_{-3}+\sum_j\partial_{\xi_j}p_3D_{x_j}q_{-3})+\cdots,
\end{align}
by (4.7), we have

\begin{equation}
q_{-3}=p_3^{-1};~q_{-4}=-p_3^{-1}[p_2p_3^{-1}+\sum_j\partial_{\xi_j}p_3D_{x_j}(p_3^{-1})].
\end{equation}
By Lemma 4.1, we have some symbols of operators.
\begin{lem} The following identities hold:
\begin{align}
\sigma_{-3}(({\widetilde{D}}^*{\widetilde{D}}{\widetilde{D}}^*)^{-1})&=\frac{ic(\xi)}{|\xi|^{4}};\nonumber\\
\sigma_{-4}(({\widetilde{D}}^*{\widetilde{D}}{\widetilde{D}}^*)^{-1})&=
\frac{c(\xi)\sigma_{2}({\widetilde{D}}^*{\widetilde{D}}{\widetilde{D}}^*)c(\xi)}{|\xi|^8}
+\frac{ic(\xi)}{|\xi|^8}\Big(|\xi|^4c(dx_n)\partial_{x_n}c(\xi')
-2h'(0)c(dx_n)c(\xi)\nonumber\\
&+2\xi_{n}c(\xi)\partial_{x_n}c(\xi')+4\xi_{n}h'(0)\Big).
\end{align}
\end{lem}
When $n=6$, then ${\rm tr}_{\wedge ^*T^*M}[\texttt{id}]=8$, where ${\rm tr}$ as shorthand of ${\rm trace}$.
Since the sum is taken over $r+\ell-k-j-|\alpha|-1=-6, \ r\leq-1, \ell\leq -3$, then we have the sum of the following five cases:

~\\
\noindent  {\bf case (a)~(I)}~$r=-1, l=-3, j=k=0, |\alpha|=1$.\\
By (4.2), we get
 \begin{equation}
\Psi_1=-\int_{|\xi'|=1}\int^{+\infty}_{-\infty}\sum_{|\alpha|=1}{\rm trace}
\Big[\partial^{\alpha}_{\xi'}\pi^{+}_{\xi_{n}}\sigma_{-1}({\widetilde{D}}^{-1})
      \times\partial^{\alpha}_{x'}\partial_{\xi_{n}}\sigma_{-3}(({\widetilde{D}}^{*}{\widetilde{D}}
      {\widetilde{D}}^{*})^{-1})\Big](x_0)d\xi_n\sigma(\xi')dx'.
\end{equation}
By Lemma 4.2, for $i<n$, we have
 \begin{equation}
 \partial_{x_{i}}\sigma_{-3}(({\widetilde{D}}^{*}{\widetilde{D}}{\widetilde{D}}^{*})^{-1})(x_0)=
      \partial_{x_{i}}\Big[\frac{ic(\xi)}{|\xi|^{4}}\Big](x_{0})
      =i\partial_{x_{i}}[c(\xi)]|\xi|^{-4}(x_{0})-2ic(\xi)\partial_{x_{i}}[|\xi|^{2}]|\xi|^{-6}(x_{0})=0,
\end{equation}
 so $\Psi_1=0$.
~\\

\noindent  {\bf case (a)~(II)}~$r=-1, l=-3, |\alpha|=k=0, j=1$.\\
By (4.2), we have
  \begin{equation}
\Psi_2=-\frac{1}{2}\int_{|\xi'|=1}\int^{+\infty}_{-\infty} {\rm
trace} \Big[\partial_{x_{n}}\pi^{+}_{\xi_{n}}\sigma_{-1}({\widetilde{D}}^{-1})
      \times\partial^{2}_{\xi_{n}}\sigma_{-3}(({\widetilde{D}}^{*}{\widetilde{D}}{\widetilde{D}}^{*})^{-1})\Big](x_0)d\xi_n\sigma(\xi')dx'.
\end{equation}
By computations, we have\\
\begin{equation}
\partial^{2}_{\xi_{n}}\sigma_{-3}(({\widetilde{D}}^{*}{\widetilde{D}}{\widetilde{D}}^{*})^{-1})=i\bigg[\frac{(20\xi^{2}_{n}-4)c(\xi')+
12(\xi^{3}_{n}-\xi_{n})c(dx_{n})}{(1+\xi_{n}^{2})^{4}}\bigg].
\end{equation}
Since $n=6$, ${\rm tr}[-\texttt{id}]=-8$. By the relation of the Clifford action and ${\rm tr}AB={\rm tr}BA$,  then
\begin{eqnarray}
&&{\rm tr}[c(\xi')c(dx_{n})]=0; \ {\rm tr}[c(dx_{n})^{2}]=-8;\
{\rm tr}[c(\xi')^{2}](x_{0})|_{|\xi'|=1}=-8;\nonumber\\
&&{\rm tr}[\partial_{x_{n}}[c(\xi')]c(\texttt{d}x_{n})]=0; \
{\rm tr}[\partial_{x_{n}}c(\xi')c(\xi')](x_{0})|_{|\xi'|=1}=-4h'(0).
\end{eqnarray}
By (3.29), (4.13) and (4.14), we get
\begin{equation}
{\rm
trace} \Big[\partial_{x_{n}}\pi^{+}_{\xi_{n}}\sigma_{-1}({\widetilde{D}}^{-1})
      \times\partial^{2}_{\xi_{n}}\sigma_{-3}(({\widetilde{D}}^{*}{\widetilde{D}}{\widetilde{D}}^{*})^{-1})\Big](x_0)
=8 h'(0)\frac{-1-3\xi_{n}i+5\xi^{2}_{n}+3i\xi^{3}_{n}}{(\xi_{n}-i)^{6}(\xi_{n}+i)^{4}}.
\end{equation}
Then we obtain

\begin{align}
\Psi_2&=-\frac{1}{2}\int_{|\xi'|=1}\int^{+\infty}_{-\infty} h'(0)dimF\frac{-1-3\xi_{n}i+5\xi^{2}_{n}+3i\xi^{3}_{n}}{(\xi_{n}-i)^{6}(\xi_{n}+i)^{4}}d\xi_n\sigma(\xi')dx'\nonumber\\
     &=h'(0)\Omega_{4}\int_{\Gamma^{+}}\frac{4+12\xi_{n}i-20\xi^{2}_{n}-122i\xi^{3}_{n}}{(\xi_{n}-i)^{6}(\xi_{n}+i)^{4}}d\xi_{n}dx'\nonumber\\
     &=h'(0)\Omega_{4}\frac{\pi i}{5!}\Big[\frac{1+3\xi_{n}i-5\xi^{2}_{n}-3i\xi^{3}_{n}}{(\xi_{n}+i)^{4}}\Big]^{(5)}|_{\xi_{n}=i}dx'\nonumber\\
     &=-\frac{15}{16}\pi h'(0)\Omega_{4}dx',\nonumber\\
\end{align}
where ${\rm \Omega_{4}}$ is the canonical volume of $S^{4}.$\\

\noindent  {\bf case (a)~(III)}~$r=-1,l=-3,|\alpha|=j=0,k=1$.\\
By (4.2), we have
 \begin{equation}
\Psi_3=-\frac{1}{2}\int_{|\xi'|=1}\int^{+\infty}_{-\infty}{\rm trace} \Big[\partial_{\xi_{n}}\pi^{+}_{\xi_{n}}\sigma_{-1}({\widetilde{D}}^{-1})
      \times\partial_{\xi_{n}}\partial_{x_{n}}\sigma_{-3}(({\widetilde{D}}^{*}{\widetilde{D}}{\widetilde{D}}^{*})^{-1})\Big](x_0)d\xi_n\sigma(\xi')dx'.
\end{equation}
By computations, we have\\
\begin{equation}
\partial_{\xi_{n}}\partial_{x_{n}}\sigma_{-3}(({\widetilde{D}}^{*}{\widetilde{D}}{\widetilde{D}}^{*})^{-1})=-\frac{4 i\xi_{n}\partial_{x_{n}}c(\xi')(x_{0})}{(1+\xi_{n}^{2})^{3}}
      +i\frac{12h'(0)\xi_{n}c(\xi')}{(1+\xi_{n}^{2})^{4}}
      -i\frac{(2-10\xi^{2}_{n})h'(0)c(dx_{n})}{(1+\xi_{n}^{2})^{4}}.
\end{equation}
Combining (3.36) and (4.18), we have
\begin{equation}
{\rm trace} \Big[\partial_{\xi_{n}}\pi^{+}_{\xi_{n}}\sigma_{-1}({\widetilde{D}}^{-1})
      \times\partial_{\xi_{n}}\partial_{x_{n}}\sigma_{-3}(({\widetilde{D}}^{*}{\widetilde{D}}{\widetilde{D}}^{*})^{-1})\Big](x_{0})|_{|\xi'|=1}
=h'(0)\frac{8i-32\xi_{n}-8i\xi^{2}_{n}}{(\xi_{n}-i)^{5}(\xi+i)^{4}}.
\end{equation}
Then
\begin{align}
\Psi_3&=-\frac{1}{2}\int_{|\xi'|=1}\int^{+\infty}_{-\infty} h'(0)\frac{8i-32\xi_{n}-8i\xi^{2}_{n}}{(\xi_{n}-i)^{5}(\xi+i)^{4}}d\xi_n\sigma(\xi')dx'\nonumber\\
     &=-\frac{1}{2}h'(0)\Omega_{4}\int_{\Gamma^{+}}\frac{8i-32\xi_{n}-8i\xi^{2}_{n}}{(\xi_{n}-i)^{5}(\xi+i)^{4}}d\xi_{n}dx'\nonumber\\
     &=-h'(0)\Omega_{4}\frac{\pi i}{4!}\Big[\frac{8i-32\xi_{n}-8i\xi^{2}_{n}}{(\xi+i)^{4}}\Big]^{(4)}|_{\xi_{n}=i}dx'\nonumber\\
     &=\frac{25}{16}\pi h'(0)\Omega_{4}dx'.
\end{align}

\noindent  {\bf case (b)}~$r=-1,l=-4,|\alpha|=j=k=0$.\\
By (4.2), we have
 \begin{align}
\Psi_4&=-i\int_{|\xi'|=1}\int^{+\infty}_{-\infty}{\rm trace} \Big[\pi^{+}_{\xi_{n}}\sigma_{-1}({\widetilde{D}}^{-1})
      \times\partial_{\xi_{n}}\sigma_{-4}(({\widetilde{D}}^{*}{\widetilde{D}}
      {\widetilde{D}}^{*})^{-1})\Big](x_0)d\xi_n\sigma(\xi')dx'\nonumber\\
&=i\int_{|\xi'|=1}\int^{+\infty}_{-\infty}{\rm trace} [\partial_{\xi_n}\pi^+_{\xi_n}\sigma_{-1}({\widetilde{D}}^{-1})\times
\sigma_{-4}(({\widetilde{D}}^{*}{\widetilde{D}}
      {\widetilde{D}}^{*})^{-1})](x_0)d\xi_n\sigma(\xi')dx'.
\end{align}

In the normal coordinate, $g^{ij}(x_{0})=\delta^{j}_{i}$ and $\partial_{x_{j}}(g^{\alpha\beta})(x_{0})=0$, if $j<n$; $\partial_{x_{j}}(g^{\alpha\beta})(x_{0})=h'(0)\delta^{\alpha}_{\beta}$, if $j=n$.
So by  \cite{Wa3}, when $k<n$, we have $\Gamma^{n}(x_{0})=\frac{5}{2}h'(0)$, $\Gamma^{k}(x_{0})=0$, $\delta^{n}(x_{0})=0$ and $\delta^{k}=\frac{1}{4}h'(0)c(e_{k})c(e_{n})$. Then, we obtain

\begin{align}
\sigma_{-4}(({\widetilde{D}}^{*}{\widetilde{D}}{\widetilde{D}}^{*})^{-1})(x_{0})|_{|\xi'|=1}&=
\frac{c(\xi)\sigma_{2}({\widetilde{D}}^{*}{\widetilde{D}}{\widetilde{D}}^{*})
(x_{0})|_{|\xi'|=1}c(\xi)}{|\xi|^8}
-\frac{c(\xi)}{|\xi|^4}\sum_j\partial_{\xi_j}\big(c(\xi)|\xi|^2\big){D}_{x_j}\big(\frac{ic(\xi)}{|\xi|^4}\big)\nonumber\\
&=\frac{1}{|\xi|^8}c(\xi)\Big(\frac{1}{2}h'(0)c(\xi)\sum_{k<n}\xi_kc(e_k)c(e_n)-\frac{5}{2}h'(0)\xi_nc(\xi)-\frac{1}{4}h'(0)\nonumber\\
&|\xi|^2c(dx_n)
+2[|\xi|^2(-\overline{f_1}\sum_{u<v}(p_{uv}-p_{vu})c(e_u)c(e_v)+\overline{f_2})-c(\xi)\nonumber\\
&(f_1\sum_{u<v}(p_{uv}-p_{vu})c(e_u)c(e_v)+f_2)c(\xi)]+|\xi|^2(-\overline{f_1}\sum_{u<v}(p_{uv}-p_{vu})\nonumber\\
&c(e_u)c(e_v)+\overline{f_2})\Big)c(\xi)+\frac{ic(\xi)}{|\xi|^8}\Big(|\xi|^4c(dx_n)\partial_{x_n}c(\xi')
-2h'(0)c(dx_n)c(\xi)\nonumber\\
&+2\xi_{n}c(\xi)\partial_{x_n}c(\xi')+4\xi_{n}h'(0)\Big).\nonumber\\
\end{align}
By (3.32) and (4.22), we have
\begin{align}
&{\rm tr} [\partial_{\xi_n}\pi^+_{\xi_n}\sigma_{-1}({\widetilde{D}}^{-1})\times
\sigma_{-4}({\widetilde{D}}^{*}{\widetilde{D}}{\widetilde{D}}^{*})^{-1}](x_0)|_{|\xi'|=1} \nonumber\\
&=\frac{1}{2(\xi_{n}-i)^{2}(1+\xi_{n}^{2})^{4}}\big(\frac{3}{4}i+2+(3+4i)\xi_{n}+(-6+2i)\xi_{n}^{2}+3\xi_{n}^{3}+\frac{9i}{4}\xi_{n}^{4}\big)h'(0){\rm tr}
[id]\nonumber\\
&+\frac{1}{2(\xi_{n}-i)^{2}(1+\xi_{n}^{2})^{4}}\big(-1-3i\xi_{n}-2\xi_{n}^{2}-4i\xi_{n}^{3}-\xi_{n}^{4}-i\xi_{n}^{5}\big){\rm tr[c(\xi')\partial_{x_n}c(\xi')]}.\nonumber\\
\end{align}
Then by (4.22), we get
\begin{align}
\Psi_4&=
 ih'(0)\int_{|\xi'|=1}\int^{+\infty}_{-\infty}8\times\frac{\frac{3}{4}i+2+(3+4i)\xi_{n}+(-6+2i)\xi_{n}^{2}+3\xi_{n}^{3}+\frac{9i}{4}\xi_{n}^{4}}{2(\xi_n-i)^5(\xi_n+i)^4}d\xi_n\sigma(\xi')dx'\nonumber\\ &+ih'(0)\int_{|\xi'|=1}\int^{+\infty}_{-\infty}4\times\frac{1+3i\xi_{n}
 +2\xi_{n}^{2}+4i\xi_{n}^{3}+\xi_{n}^{4}+i\xi_{n}^{5}}{2(\xi_{n}-i)^{2}
 (1+\xi_{n}^{2})^{4}}d\xi_n\sigma(\xi')dx'\nonumber\\
&=(-\frac{41}{64}i-\frac{195}{64})\pi h'(0)\Omega_4dx'.\nonumber\\
\end{align}

\noindent {\bf  case (c)}~$r=-2,l=-3,|\alpha|=j=k=0$.\\
By (4.2), we have

\begin{equation}
\Psi_5=-i\int_{|\xi'|=1}\int^{+\infty}_{-\infty}{\rm trace} \Big[\pi^{+}_{\xi_{n}}\sigma_{-2}({\widetilde{D}}^{-1})
      \times\partial_{\xi_{n}}\sigma_{-3}(({\widetilde{D}}^{*}{\widetilde{D}}{\widetilde{D}}^{*})^{-1})\Big](x_0)d\xi_n\sigma(\xi')dx'.
\end{equation}

By Lemma 4.1 and Lemma 4.2, we have
\begin{align}
\sigma_{-2}({\widetilde{D}}^{-1})(x_0)&=\frac{c(\xi)\sigma_{0}({\widetilde{D}})c(\xi)}{|\xi|^4}(x_0)+\frac{c(\xi)}{|\xi|^6}\sum_jc(dx_j)
\Big[\partial_{x_j}(c(\xi))|\xi|^2-c(\xi)\partial_{x_j}(|\xi|^2)\Big](x_0),
\end{align}
where
\begin{align}
\sigma_0({\widetilde{D}})=
-\frac{1}{4}\sum_{i,s,t}\omega_{s,t}(e_i)c(e_i)c(e_s)c(e_t)+(f_1\sum_{u<v}(p_{uv}-p_{vu})c(e_u)c(e_v)+f_2).
\end{align}

On the other hand,
\begin{align}
\partial_{\xi_{n}}\sigma_{-3}(({\widetilde{D}}^{*}{\widetilde{D}}{\widetilde{D}}^{*})^{-1})=\frac{-4 i \xi_{n}c(\xi')}{(1+\xi_{n}^{2})^{3}}+\frac{i(1- 3\xi_{n}^{2})c(\texttt{d}x_{n})}
{(1+\xi_{n}^{2})^{3}}.
\end{align}
By (4.26), (3.5) and (3.6), we have
\begin{align}
\pi^{+}_{\xi_{n}}\Big(\sigma_{-2}({\widetilde{D}}^{-1})\Big)(x_{_{0}})|_{|\xi'|=1}
&=\pi^{+}_{\xi_{n}}\Big[\frac{c(\xi)\sigma_{0}({\widetilde{D}})(x_{0})c(\xi)
+c(\xi)c(dx_{n})\partial_{x_{n}}[c(\xi')](x_{0})}{(1+\xi^{2}_{n})^{2}}\Big]\nonumber\\
&-h'(0)\pi^{+}_{\xi_{n}}\Big[\frac{c(\xi)c(dx_{n})c(\xi)}{(1+\xi^{2}_{n})^{3}}\Big].
\end{align}
We denote
 \begin{align}
\sigma_{0}({\widetilde{D}})(x_0)|_{\xi_n=i}=Q(x_0)
+(f_1\sum_{u<v}(p_{uv}-p_{vu})c(e_u)c(e_v)+f_2).
\end{align}
Then, we obtain
\begin{align}
\pi^{+}_{\xi_{n}}\Big(\sigma_{-2}({\widetilde{D}}^{-1})\Big)(x_{_{0}})|_{|\xi'|=1}
&=\pi^+_{\xi_n}\Big[\frac{c(\xi)Q(x_0)c(\xi)+c(\xi)c(dx_n)
\partial_{x_n}[c(\xi')](x_0)}{(1+\xi_n^2)^2}-h'(0)\frac{c(\xi)c(dx_n)c(\xi)}{(1+\xi_n^{2})^3}\Big]\nonumber\\
&+\pi^+_{\xi_n}\Big[\frac{c(\xi)[(f_1\sum_{u<v}(p_{uv}-p_{vu})c(e_u)c(e_v)+f_2)]c(\xi)(x_0)}{(1+\xi_n^2)^2}\Big].\nonumber\\
\end{align}

Furthermore,
\begin{align}
&\pi^+_{\xi_n}\Big[\frac{c(\xi)[(f_1\sum_{u<v}(p_{uv}-p_{vu})c(e_u)c(e_v)+f_2)](x_0)c(\xi)}
{(1+\xi_n^2)^2}\Big]\nonumber\\
&=\pi^+_{\xi_n}\Big[\frac{c(\xi')[(f_1\sum_{u<v}(p_{uv}-p_{vu})c(e_u)c(e_v)+f_2)](x_0)c(\xi')}
{(1+\xi_n^2)^2}\Big]\nonumber\\
&
+\pi^+_{\xi_n}\Big[ \frac{\xi_nc(\xi')[(f_1\sum_{u<v}(p_{uv}-p_{vu})c(e_u)c(e_v)+f_2)]
(x_0)c(dx_{n})}{(1+\xi_n^2)^2}\Big]\nonumber\\
&+\pi^+_{\xi_n}\Big[\frac{\xi_nc(dx_{n})[(f_1\sum_{u<v}(p_{uv}-p_{vu})c(e_u)c(e_v)+f_2)](x_0)c(\xi')}{(1+\xi_n^2)^2}\Big]\nonumber\\
&
+\pi^+_{\xi_n}\Big[\frac{\xi_n^{2}c(dx_{n})[(f_1\sum_{u<v}(p_{uv}-p_{vu})c(e_u)c(e_v)+f_2)](x_0)c(dx_{n})}{(1+\xi_n^2)^2}\Big]\nonumber\\
&=-\frac{c(\xi')[(f_1\sum_{u<v}(p_{uv}-p_{vu})c(e_u)c(e_v)+f_2)](x_0)c(\xi')(2+i\xi_{n})}{4(\xi_{n}-i)^{2}}\nonumber\\
&
+\frac{ic(\xi')[(f_1\sum_{u<v}(p_{uv}-p_{vu})c(e_u)c(e_v)+f_2)](x_0)c(dx_{n})}{4(\xi_{n}-i)^{2}}\nonumber\\
&+\frac{ic(dx_{n})[(f_1\sum_{u<v}(p_{uv}-p_{vu})c(e_u)c(e_v)+f_2)](x_0)c(\xi')}{4(\xi_{n}-i)^{2}}\nonumber\\
&
+\frac{-i\xi_{n}c(dx_{n})[(f_1\sum_{u<v}(p_{uv}-p_{vu})c(e_u)c(e_v)+f_2)](x_0)c(dx_{n})}{4(\xi_{n}-i)^{2}},\nonumber\\
\end{align}
By computations, we have
\begin{align}
\pi^+_{\xi_n}\Big[\frac{c(\xi)Q(x_0)c(\xi)+c(\xi)c(dx_n)\partial_{x_n}(c(\xi'))(x_0)}{(1+\xi_n^2)^2}\Big]-h'(0)\pi^+_{\xi_n}\Big[\frac{c(\xi)c(dx_n)c(\xi)}{(1+\xi_n)^3}\Big]:= C_1-C_2,\nonumber\\
\end{align}
where
\begin{align}
C_1&=\frac{-1}{4(\xi_n-i)^2}\big[(2+i\xi_n)c(\xi')Qc(\xi')+i\xi_nc(dx_n)Qc(dx_n) \nonumber\\
&+(2+i\xi_n)c(\xi')c(dx_n)\partial_{x_n}c(\xi')+ic(dx_n)Q^2_0c(\xi')
+ic(\xi')Qc(dx_n)-i\partial_{x_n}c(\xi')\big]\nonumber\\
&=\frac{1}{4(\xi_n-i)^2}\Big[\frac{5}{2}h'(0)c(dx_n)-\frac{5i}{2}h'(0)c(\xi')
  -(2+i\xi_n)c(\xi')c(dx_n)\partial_{\xi_n}c(\xi')+i\partial_{\xi_n}c(\xi')\Big]; \nonumber\\
C_2&=\frac{h'(0)}{2}\Big[\frac{c(dx_n)}{4i(\xi_n-i)}+\frac{c(dx_n)-ic(\xi')}{8(\xi_n-i)^2}
+\frac{3\xi_n-7i}{8(\xi_n-i)^3}\big(ic(\xi')-c(dx_n)\big)\Big].\nonumber\\
\end{align}
By (4.29) and (4.34), we have
\begin{align}
&{\rm tr }[C_2\times\partial_{\xi_n}\sigma_{-3}(({\widetilde{D}}^{*}{\widetilde{D}}{\widetilde{D}}^{*})^{-1})(x_0)]|_{|\xi'|=1}\nonumber\\
&={\rm tr }\Big\{ \frac{h'(0)}{2}\Big[\frac{c(dx_n)}{4i(\xi_n-i)}+\frac{c(dx_n)-ic(\xi')}{8(\xi_n-i)^2}
+\frac{3\xi_n-7i}{8(\xi_n-i)^3}[ic(\xi')-c(dx_n)]\Big] \nonumber\\
&\times\frac{-4i\xi_nc(\xi')+(i-3i\xi_n^{2})c(dx_n)}{(1+\xi_n^{2})^3}\Big\} \nonumber\\
&=h'(0)\frac{4i-11\xi_n-6i\xi_n^{2}+3\xi_n^{3}}{(\xi_n-i)^5(\xi_n+i)^3}.
\end{align}
Similarly, we have
\begin{align}
&{\rm tr }[C_1\times\partial_{\xi_n}\sigma_{-3}(({\widetilde{D}}^{*}{\widetilde{D}}{\widetilde{D}}^{*})^{-1})(x_0)]|_{|\xi'|=1}\nonumber\\
&={\rm tr }\Big\{ \frac{1}{4(\xi_n-i)^2}\Big[\frac{5}{2}h'(0)c(dx_n)-\frac{5i}{2}h'(0)c(\xi')
  -(2+i\xi_n)c(\xi')c(dx_n)\partial_{\xi_n}c(\xi')+i\partial_{\xi_n}c(\xi')\Big]\nonumber\\
&\times \frac{-4i\xi_nc(\xi')+(i-3i\xi_n^{2})c(dx_n)}{(1+\xi_n^{2})^3}\Big\} \nonumber\\
&=h'(0)\frac{3+12i\xi_n+3\xi_n^{2}}{(\xi_n-i)^4(\xi_n+i)^3};\nonumber\\
&{\rm tr }\bigg[\pi^+_{\xi_n}\Big(\frac{c(\xi)[(f_1\sum_{u<v}(p_{uv}-p_{vu})c(e_u)c(e_v)+f_2)]
(x_0)c(\xi)}{(1+\xi_n^2)^2}\Big)\times
\partial_{\xi_n}\sigma_{-3}(({D}^{*}{\widetilde{D}}{\widetilde{D}}^{*})^{-1})
(x_0)\bigg]\bigg|_{|\xi'|=1}\nonumber\\
&=\frac{2-8i\xi_n-6\xi_n^2}{4(\xi_n-i)^{2}(1+\xi_n^2)^{3}}{\rm tr }[(f_1\sum_{u<v}(p_{uv}-p_{vu})c(e_u)c(e_v)+f_2)(x_0)c(\xi')]\nonumber\\
&=0.\nonumber\\
\end{align}
By $\int_{|\xi'|=1}\xi_{1}\cdot\cdot\cdot\xi_{2q+1}\sigma(\xi')=0,$ we have\\
\begin{align}
\Psi_5&=
 -i h'(0)\int_{|\xi'|=1}\int^{+\infty}_{-\infty}
 \times\frac{-7i+26\xi_n+15i\xi_n^{2}}{(\xi_n-i)^5(\xi_n+i)^3}d\xi_n\sigma(\xi')dx' \nonumber\\
&=-i h'(0)\times\frac{2 \pi i}{4!}\Big[\frac{-7i+26\xi_n+15i\xi_n^{2}}{(\xi_n+i)^3}
     \Big]^{(5)}|_{\xi_n=i}\Omega_4dx'\nonumber\\
&=\frac{55}{16}\pi h'(0)\Omega_4dx'.
\end{align}

Now $\Psi$ is the sum of the cases (a), (b) and (c), then
\begin{equation}
\Psi=\Psi_1+\Psi_2+\Psi_3+\Psi_4+\Psi_5=(\frac{65}{64}-\frac{41}{64}i)\pi h'(0)\Omega_4dx'.
\end{equation}

\begin{thm}
Let $M$ be a $6$-dimensional
compact oriented spin manifold with the boundary $\partial M$ and the metric
$g^M$ as above, ${\widetilde{D}}$ and ${\widetilde{D}}^*$ be the Dirac-Witten operators on $\widetilde{M}$, then
\begin{align}
&\widetilde{{\rm Wres}}[\pi^+{\widetilde{D}}^{-1}\circ\pi^+({\widetilde{D}}^{*}{\widetilde{D}}
      {\widetilde{D}}^{*})^{-1}]\nonumber\\
&=128\pi^3\int_{M}\bigg(
-\frac{2}{3}s-4(f_1^2+\overline{f_1}^2-4f_1\overline{f_1})\sum_{u<v}(p_{uv}-p_{vu})^2-32f_2\overline{f_2}+12f_2^2+12\overline{f_2}^2\bigg)d{\rm Vol_{M}}\nonumber\\
&+\int_{\partial M}(\frac{65}{64}-\frac{41}{64}i)\pi h'(0)\Omega_4dx'.\nonumber\\
\end{align}
where $s$ is the scalar curvature.
\end{thm}

Next, we prove the Kastler-Kalau-Walze type theorem for $6$-dimensional manifold with boundary associated to ${\widetilde{D}}^{3}$. From \cite{Wa5}, we know that

\begin{equation}
\widetilde{{\rm Wres}}[\pi^+{\widetilde{D}}^{-1}\circ\pi^+{\widetilde{D}}^{-3}
      ]=\int_M\int_{|\xi|=1}{\rm
trace}_{\wedge^*T^*M\bigotimes\mathbb{C}}[\sigma_{-4}({\widetilde{D}}^{-4})]\sigma(\xi)dx+\int_{\partial M}\overline{\Psi},
\end{equation}
where $\widetilde{{\rm Wres}}$ denote noncommutative residue on minifolds with boundary,
\begin{align}
\overline{\Psi} &=\int_{|\xi'|=1}\int^{+\infty}_{-\infty}\sum^{\infty}_{j, k=0}\sum\frac{(-i)^{|\alpha|+j+k+1}}{\alpha!(j+k+1)!}
\times {\rm trace}_{{\wedge^*T^*M\bigotimes\mathbb{C}}}[\partial^j_{x_n}\partial^\alpha_{\xi'}\partial^k_{\xi_n}\sigma^+_{r}({\widetilde{D}}^{-1})(x',0,\xi',\xi_n)
\nonumber\\
&\times\partial^\alpha_{x'}\partial^{j+1}_{\xi_n}\partial^k_{x_n}\sigma_{l}
({\widetilde{D}}^{-3})(x',0,\xi',\xi_n)]d\xi_n\sigma(\xi')dx',\nonumber\\
\end{align}
and the sum is taken over $r+\ell-k-j-|\alpha|-1=-6, \ r\leq-1, \ell\leq -3$.

By Theorem 2.2, we compute the interior term of (4.41), then
\begin{align}
&\int_M\int_{|\xi|=1}{\rm
trace}_{\wedge^*T^*M\bigotimes\mathbb{C}}[\sigma_{-4}({\widetilde{D}}^{-4})]\sigma(\xi)dx\nonumber\\
&=128\pi^3\int_{M}\bigg[-\frac{2}{3}s
-24f_1^2\sum_{u<v}(p_{uv}-p_{vu})^2+40f_2^2\bigg]d{\rm Vol_{M}}\nonumber\\.
\end{align}

So we only need to compute $\int_{\partial M} \overline{\Psi}$. Let us now turn to compute the specification of
${\widetilde{D}}^3$.
\begin{align}
{\widetilde{D}}^3
&=\sum^{n}_{i=1}c(e_{i})\langle e_{i},dx_{l}\rangle(-g^{ij}\partial_{l}\partial_{i}\partial_{j})
+\sum^{n}_{i=1}c(e_{i})\langle e_{i},dx_{l}\rangle \bigg\{-(\partial_{l}g^{ij})\partial_{i}\partial_{j}-g^{ij}\bigg(4\sigma_{i}
\partial_{j}2\Gamma^{k}_{ij}\partial_{k}\bigg)\nonumber\\
&\partial_{l}\bigg\}+\sum^{n}_{i=1}c(e_{i})\langle e_{i},dx_{l}\rangle \bigg\{-2(\partial_{l}g^{ij})\sigma_{i}\partial_{j}+g^{ij}
(\partial_{l}\Gamma^{k}_{ij})\partial_{k}-2g^{ij}[(\partial_{l}\sigma_{i})
+(\partial_{l}g^{ij})\Gamma^{k}_{ij}\partial_{k}\nonumber\\
&+\sum_{j,k}\Big[\partial_{l}\Big((f_1\sum_{u<v}(p_{uv}-p_{vu})c(e_u)c(e_v)+f_2)c(e_{j})+c(e_{j})(f_1\sum_{u<v}(p_{uv}-p_{vu})c(e_u)c(e_v)\nonumber\\
&+f_2)\Big)\Big]\langle e_{j},dx^{k}\rangle\partial_{k}+\sum_{j,k}\Big((f_1\sum_{u<v}(p_{uv}-p_{vu})c(e_u)c(e_v)+f_2)c(e_{j})+c(e_{j})(f_1\sum_{u<v}(p_{uv}\nonumber\\
&-p_{vu})c(e_u)c(e_v)+f_2)\Big)\Big[\partial_{l}\langle e_{j},dx^{k}\rangle\Big]\partial_{k} \bigg\}+\sum^{n}_{i=1}c(e_{i})\langle e_{i},dx_{l}\rangle\partial_{l}\bigg\{-g^{ij}\Big[(\partial_{i}\sigma_{j})+\sigma_{i}\sigma_{j}\nonumber\\
&-\Gamma_{i,j}^{k}\sigma_{k}+\sum_{i,j}g^{i,j}\Big[(f_1\sum_{u<v}(p_{uv}-p_{vu})c(e_u)c(e_v)+f_2)c(\partial_{i})\sigma_{i}+c(\partial_{i})\partial_{i}(f_1\sum_{u<v}(p_{uv}-p_{vu})\nonumber\\
&c(e_u)c(e_v)+f_2)+c(\partial_{i})\sigma_{i}(f_1\sum_{u<v}(p_{uv}-p_{vu})c(e_u)c(e_v)+f_2)\Big]+\frac{1}{4}s-[(f_1\sum_{u<v}(p_{uv}-p_{vu})\nonumber\\
&c(e_u)c(e_v)+f_2)]^2\bigg\}+\Big[\sigma_{i}+(f_1\sum_{u<v}(p_{uv}-p_{vu})c(e_u)c(e_v)+f_2)\Big](-g^{ij}\partial_{i}\partial_{j})+\sum^{n}_{i=1}c(e_{i})\nonumber\\
&\langle e_{i},dx_{l}\rangle\bigg\{2\sum_{j,k}\Big[(f_1\sum_{u<v}(p_{uv}-p_{vu})c(e_u)c(e_v)+f_2)c(e_{j})+c(e_{j})(f_1\sum_{u<v}(p_{uv}-p_{vu})c(e_u)\nonumber\\
&c(e_v)+f_2)\Big]\times\langle e_{i},dx_{k}\rangle\bigg\}_{l}\partial_{k}
+\Big[\sigma_{i}+(f_1\sum_{u<v}(p_{uv}-p_{vu})c(e_u)c(e_v)+f_2)\Big]
\bigg\{-\sum_{i,j}g^{i,j}\nonumber\\
&\Big[2\sigma_{i}\partial_{j}
-\Gamma_{i,j}^{k}\partial_{k}+(\partial_{i}\sigma_{j})
+\sigma_{i}\sigma_{j} -\Gamma_{i,j}^{k}\sigma_{k}\Big]+\sum_{i,j}g^{i,j}\Big[c(\partial_{i})(f_1\sum_{u<v}(p_{uv}-p_{vu})c(e_u)c(e_v)\nonumber\\
&+f_2)
+(f_1\sum_{u<v}(p_{uv}-p_{vu})c(e_u)c(e_v)+f_2)c(\partial_{i})\Big]\partial_{j}+\sum_{i,j}g^{i,j}\Big[(f_1\sum_{u<v}(p_{uv}-p_{vu})c(e_u)c(e_v)\nonumber\\
&+f_2)c(\partial_{i})\sigma_{i}+c(\partial_{i})\partial_{i}(f_1\sum_{u<v}(p_{uv}-p_{vu})c(e_u)c(e_v)+f_2)+c(\partial_{i})\sigma_{i}(f_1\sum_{u<v}(p_{uv}-p_{vu})c(e_u)\nonumber\\
&c(e_v)+f_2)+c(\partial_{i})\partial_{i}(f_1\sum_{u<v}(p_{uv}-p_{vu})c(e_u)c(e_v)+f_2)+c(\partial_{i})\sigma_{i}(f_1\sum_{u<v}(p_{uv}-p_{vu})c(e_u)\nonumber\\
&c(e_v)+f_2)\Big]+\frac{1}{4}s-[(f_1\sum_{u<v}(p_{uv}-p_{vu})c(e_u)c(e_v)+f_2)]^2\bigg\}.\nonumber\\
\end{align}

Then, we obtain
\begin{lem} The following identities hold:
\begin{align}
\sigma_2({\widetilde{D}}^3)&=\sum_{i,j,l}c(dx_{l})\partial_{l}(g^{i,j})\xi_{i}\xi_{j} +c(\xi)(4\sigma^k-2\Gamma^k)\xi_{k}-2[c(\xi)(f_1\sum_{u<v}(p_{uv}-p_{vu})\nonumber\\
&c(e_u)c(e_v)+f_2)c(\xi)-|\xi|^2(f_1\sum_{u<v}(p_{uv}-p_{vu})c(e_u)c(e_v)+f_2)]\nonumber\\
&-\frac{1}{4}|\xi|^2\sum_{s,t,l}\omega_{s,t}
(e_l)c(e_s)c(e_t)];\nonumber\\
\sigma_{3}({\widetilde{D}}^3)&=ic(\xi)|\xi|^{2}.\nonumber\\
\end{align}
\end{lem}

Write
\begin{eqnarray}
\sigma({\widetilde{D}}^3)=p_3+p_2+p_1+p_0;
~\sigma({\widetilde{D}}^{-3})=\sum^{\infty}_{j=3}q_{-j}.
\end{eqnarray}

By the composition formula of pseudodifferential operators, we have

\begin{align}
1=\sigma({\widetilde{D}}^3\circ {\widetilde{D}}^{-3})&=
\sum_{\alpha}\frac{1}{\alpha!}\partial^{\alpha}_{\xi}
[\sigma({\widetilde{D}}^3)]{\widetilde{D}}^{\alpha}_{x}
[\sigma({\widetilde{D}}^{-3})] \nonumber\\
&=(p_3+p_2+p_1+p_0)(q_{-3}+q_{-4}+q_{-5}+\cdots) \nonumber\\
&+\sum_j(\partial_{\xi_j}p_3+\partial_{\xi_j}p_2+\partial_{\xi_j}p_1+\partial_{\xi_j}p_0)
(D_{x_j}q_{-3}+D_{x_j}q_{-4}+D_{x_j}q_{-5}+\cdots) \nonumber\\
&=p_3q_{-3}+(p_3q_{-4}+p_2q_{-3}+\sum_j\partial_{\xi_j}p_3D_{x_j}q_{-3})+\cdots,
\end{align}
by (4.46), we have

\begin{equation}
q_{-3}=p_3^{-1};~q_{-4}=-p_3^{-1}[p_2p_3^{-1}+\sum_j\partial_{\xi_j}p_3D_{x_j}(p_3^{-1})].
\end{equation}
By (4.43)-(4.47), we have some symbols of operators.
\begin{lem} The following identities hold:
\begin{align}
\sigma_{-3}({\widetilde{D}}^{-3})&=\frac{ic(\xi)}{|\xi|^{4}};\nonumber\\
\sigma_{-4}({\widetilde{D}}^{-3})&=
\frac{c(\xi)\sigma_{2}({\widetilde{D}}^{3})c(\xi)}{|\xi|^8}
+\frac{ic(\xi)}{|\xi|^8}\Big(|\xi|^4c(dx_n)\partial_{x_n}c(\xi')
-2h'(0)c(dx_n)c(\xi)\nonumber\\
&+2\xi_{n}c(\xi)\partial_{x_n}c(\xi')+4\xi_{n}h'(0)\Big).
\end{align}
\end{lem}

When $n=6$, then ${\rm tr}_{\wedge^*T^*M}[\texttt{id}]=8$, where ${\rm tr}$ as shorthand of ${\rm trace}$.
Since the sum is taken over $r+\ell-k-j-|\alpha|-1=-6, \ r\leq-1, \ell\leq -3$, then we have the following five cases:

~\\
\noindent  {\bf case (a)~(I)}~$r=-1, l=-3, j=k=0, |\alpha|=1$.\\
By (4.41), we get
 \begin{equation}
\overline{\Psi}_1=-\int_{|\xi'|=1}\int^{+\infty}_{-\infty}\sum_{|\alpha|=1}{\rm trace}
\Big[\partial^{\alpha}_{\xi'}\pi^{+}_{\xi_{n}}\sigma_{-1}({\widetilde{D}}^{-1})
      \times\partial^{\alpha}_{x'}\partial_{\xi_{n}}\sigma_{-3}({\widetilde{D}}^{-3})\Big](x_0)d\xi_n\sigma(\xi')dx'.
\end{equation}

\noindent  {\bf case (a)~(II)}~$r=-1, l=-3, |\alpha|=k=0, j=1$.\\
By (4.41), we have
  \begin{equation}
\overline{\Psi}_2=-\frac{1}{2}\int_{|\xi'|=1}\int^{+\infty}_{-\infty} {\rm
trace} \Big[\partial_{x_{n}}\pi^{+}_{\xi_{n}}\sigma_{-1}({\widetilde{D}}^{-1})
      \times\partial^{2}_{\xi_{n}}\sigma_{-3}({\widetilde{D}}^{-3})\Big](x_0)d\xi_n\sigma(\xi')dx'.
\end{equation}

\noindent  {\bf case (a)~(III)}~$r=-1,l=-3,|\alpha|=j=0,k=1$.\\
By (4.41), we have
 \begin{equation}
\overline{\Psi}_3=-\frac{1}{2}\int_{|\xi'|=1}\int^{+\infty}_{-\infty}{\rm trace} \Big[\partial_{\xi_{n}}\pi^{+}_{\xi_{n}}\sigma_{-1}({\widetilde{D}}^{-1})
      \times\partial_{\xi_{n}}\partial_{x_{n}}\sigma_{-3}({\widetilde{D}}^{-3})\Big](x_0)d\xi_n\sigma(\xi')dx'.
\end{equation}\\
By Lemma 4.2 and Lemma 4.5, we have $\sigma_{-3}(({\widetilde{D}}^*{\widetilde{D}}{\widetilde{D}}^*)^{-1})=\sigma_{-3}({\widetilde{D}}^{-3})$
, by (4.10)-(4.20), we obtain
$$\overline{\Psi}_1+\overline{\Psi}_2+\overline{\Psi}_3=\frac{5}{8}\pi h'(0)\Omega_{4}dx',$$
 where ${\rm \Omega_{4}}$ is the canonical volume of $S^{4}.$\\

\noindent  {\bf case (b)}~$r=-1,l=-4,|\alpha|=j=k=0$.\\
By (4.41), we have
 \begin{align}
\overline{\Psi}_4&=-i\int_{|\xi'|=1}\int^{+\infty}_{-\infty}{\rm trace} \Big[\pi^{+}_{\xi_{n}}\sigma_{-1}({\widetilde{D}}^{-1})
      \times\partial_{\xi_{n}}\sigma_{-4}({\widetilde{D}}^{-3})\Big](x_0)d\xi_n\sigma(\xi')dx'\nonumber\\
&=i\int_{|\xi'|=1}\int^{+\infty}_{-\infty}{\rm trace} [\partial_{\xi_n}\pi^+_{\xi_n}\sigma_{-1}({\widetilde{D}}^{-1})\times
\sigma_{-4}({\widetilde{D}}^{-3})](x_0)d\xi_n\sigma(\xi')dx'.
\end{align}

In the normal coordinate, $g^{ij}(x_{0})=\delta^{j}_{i}$ and $\partial_{x_{j}}(g^{\alpha\beta})(x_{0})=0$, if $j<n$; $\partial_{x_{j}}(g^{\alpha\beta})(x_{0})=h'(0)\delta^{\alpha}_{\beta}$, if $j=n$.
So by  \cite{Wa3}, when $k<n$, we have $\Gamma^{n}(x_{0})=\frac{5}{2}h'(0)$, $\Gamma^{k}(x_{0})=0$, $\delta^{n}(x_{0})=0$ and $\delta^{k}=\frac{1}{4}h'(0)c(e_{k})c(e_{n})$. Then, we obtain

\begin{align}
\sigma_{-4}({\widetilde{D}}^{-3})(x_{0})|_{|\xi'|=1}&=
\frac{c(\xi)\sigma_{2}({\widetilde{D}}^{3})
(x_{0})|_{|\xi'|=1}c(\xi)}{|\xi|^8}
-\frac{c(\xi)}{|\xi|^4}\sum_j\partial_{\xi_j}\big(c(\xi)|\xi|^2\big)
{D}_{x_j}\big(\frac{ic(\xi)}{|\xi|^4}\big)\nonumber\\
&=\frac{1}{|\xi|^8}c(\xi)\Big(\frac{1}{2}h'(0)c(\xi)\sum_{k<n}\xi_k
c(e_k)c(e_n)-\frac{5}{2}h'(0)\xi_nc(\xi)-\frac{1}{4}h'(0)|\xi|^2\nonumber\\
&c(dx_n)
-2[c(\xi)(f_1\sum_{u<v}(p_{uv}-p_{vu})c(e_u)c(e_v)+f_2)c(\xi)+|\xi|^2(-\overline{f_1}\sum_{u<v}\nonumber\\
&(p_{uv}-p_{vu})c(e_u)c(e_v)+\overline{f_2})]+|\xi|^2(f_1\sum_{u<v}(p_{uv}-p_{vu})c(e_u)c(e_v)+f_2)\Big)\nonumber\\
&c(\xi)+\frac{ic(\xi)}{|\xi|^8}\Big(|\xi|^4c(dx_n)\partial_{x_n}c(\xi')
-2h'(0)c(dx_n)c(\xi)+2\xi_{n}c(\xi)\partial_{x_n}c(\xi')\nonumber\\
&+4\xi_{n}h'(0)\Big).\nonumber\\
\end{align}
By (3.29) and (4.53), we have
\begin{align}
&{\rm tr} [\partial_{\xi_n}\pi^+_{\xi_n}\sigma_{-1}({\widetilde{D}}^{-1})\times
\sigma_{-4}({\widetilde{D}}^{-3}) ](x_0)|_{|\xi'|=1} \nonumber\\
&=\frac{1}{2(\xi_{n}-i)^{2}(1+\xi_{n}^{2})^{4}}\big(\frac{3}{4}i+2+(3+4i)\xi_{n}+(-6+2i)\xi_{n}^{2}+3\xi_{n}^{3}+\frac{9i}{4}\xi_{n}^{4}\big)h'(0){\rm tr}
[id]\nonumber\\
&+\frac{1}{2(\xi_{n}-i)^{2}(1+\xi_{n}^{2})^{4}}\big(-1-3i\xi_{n}-2\xi_{n}^{2}-4i\xi_{n}^{3}-\xi_{n}^{4}-i\xi_{n}^{5}\big){\rm tr[c(\xi')\partial_{x_n}c(\xi')]}.\nonumber\\
\end{align}
By (4.54), we have
\begin{align}
\overline{\Psi}_4&=
 ih'(0)\int_{|\xi'|=1}\int^{+\infty}_{-\infty}8\times\frac{\frac{3}{4}i+2+(3+4i)\xi_{n}+(-6+2i)\xi_{n}^{2}+3\xi_{n}^{3}+\frac{9i}{4}\xi_{n}^{4}}{2(\xi_n-i)^5(\xi_n+i)^4}d\xi_n\sigma(\xi')dx'\nonumber\\ &+ih'(0)\int_{|\xi'|=1}\int^{+\infty}_{-\infty}4\times\frac{1+3i\xi_{n}
 +2\xi_{n}^{2}+4i\xi_{n}^{3}+\xi_{n}^{4}+i\xi_{n}^{5}}{2(\xi_{n}-i)^{2}
 (1+\xi_{n}^{2})^{4}}d\xi_n\sigma(\xi')dx'\nonumber\\
&=(-\frac{41}{64}i-\frac{195}{64})\pi h'(0)\Omega_4dx'.\nonumber\\
\end{align}

\noindent {\bf  case (c)}~$r=-2,l=-3,|\alpha|=j=k=0$.\\
By (4.41), we have

\begin{equation}
\overline{\Psi}_5=-i\int_{|\xi'|=1}\int^{+\infty}_{-\infty}{\rm trace} \Big[\pi^{+}_{\xi_{n}}\sigma_{-2}({\widetilde{D}}^{-1})
      \times\partial_{\xi_{n}}\sigma_{-3}({\widetilde{D}}^{-3})\Big](x_0)d\xi_n\sigma(\xi')dx'.
\end{equation}\\
By Lemma 4.2 and Lemma 4.5, we have $\sigma_{-3}(({\widetilde{D}}^*{\widetilde{D}}{\widetilde{D}}^*)^{-1})=\sigma_{-3}({\widetilde{D}}^{-3})$
, by (4.25)-(4.37), we obtain

$$
\overline{\Psi}_5=\frac{55}{16}\pi h'(0)\Omega_4dx'.
$$

Now $\overline{\Psi}$ is the sum of the cases (a), (b) and (c), then
\begin{equation}
\overline{\Psi}=\overline{\Psi}_1+\overline{\Psi}_2+\overline{\Psi}_3+\overline{\Psi}_4+\overline{\Psi}_5=(\frac{65}{64}-\frac{41}{64}i)\pi h'(0)\Omega_4dx'.
\end{equation}

\begin{thm}
Let $M$ be a $6$-dimensional
compact oriented spin manifold with the boundary $\partial M$ and the metric
$g^M$ as above, ${\widetilde{D}}$ be the Dirac-Witten operator on $\widetilde{M}$, then
\begin{align}
&\widetilde{{\rm Wres}}[\pi^+{\widetilde{D}}^{-1}\circ\pi^+(
      {\widetilde{D}}^{-3})]\nonumber\\
&=128\pi^3\int_{M}\bigg[-\frac{2}{3}s
-24f_1^2\sum_{u<v}(p_{uv}-p_{vu})^2+40f_2^2\bigg]d{\rm Vol_{M}}+\int_{\partial M}(\frac{65}{64}-\frac{41}{64}i)\pi h'(0)\Omega_4dx'\nonumber\\.
\end{align}
where $s$ is the scalar curvature.
\end{thm}

\section*{Acknowledgements}
The author was supported in part by  NSFC No.11771070. The author thanks the referee for his (or her) careful reading and helpful comments.

\end{document}